\pdfoutput=1
\documentclass[11pt,a4paper]{article}
\usepackage[T1]{fontenc}
\usepackage[utf8]{inputenc}
\usepackage{lmodern}
\usepackage{amsmath,amssymb,amsthm}
\usepackage{graphicx}
\usepackage{booktabs}
\usepackage{microtype}
\usepackage[margin=2.6cm]{geometry}
\usepackage{xcolor}
\usepackage{pgfplots}
\usepackage{listings}
\usepackage[colorlinks=true,linkcolor=black,citecolor=black,urlcolor=blue!50!black]{hyperref}
\usepackage{caption}
\usepackage[section,above,below]{placeins}        

\newtheorem{theorem}{Theorem}
\newtheorem{proposition}[theorem]{Proposition}
\newtheorem{corollary}[theorem]{Corollary}
\newtheorem{lemma}[theorem]{Lemma}
\theoremstyle{definition}
\newtheorem{definition}[theorem]{Definition}
\newtheorem{result}[theorem]{Computational result}

\theoremstyle{remark}

\newcommand{\R}{\mathbb{R}}
\newcommand{\Z}{\mathbb{Z}}
\newcommand{\Q}{\mathbb{Q}}
\newcommand{\Zt}{\mathbb{Z}[\tau]}
\newcommand{\M}{\mathcal{M}}
\newcommand{\Pm}{\mathcal{P}}
\newcommand{\B}{B}
\newcommand{\Dz}{D_0}
\newcommand{\vol}{\operatorname{vol}}
\newcommand{\sgn}{\operatorname{sgn}}
\newcommand{\GL}{\operatorname{GL}}
\newcommand{\Om}{\Omega}
\newcommand{\dimH}{\dim_{\mathrm H}}
\newcommand{\dimB}{\overline{\dim}_{\mathrm B}}
\newcommand{\betaval}{3.192582}
\newcommand{\bka}{31}
\newcommand{\bkb}{32}
\newcommand{\bkc}{212}
\newcommand{\bkd}{392}
\newcommand{\bke}{1\,832}
\newcommand{\bkf}{4\,532}
\newcommand{\bkg}{17\,312}
\newcommand{\bkh}{48\,992}
\newcommand{\bki}{170\,132}
\newcommand{\bkj}{513\,032}
\newcommand{\bkk}{1\,703\,912}
\newcommand{\bkl}{5\,295\,092}
\newcommand{\bkm}{17\,222\,432}
\newcommand{\boxa}{540}
\newcommand{\boxb}{2\,510}
\newcommand{\boxc}{12\,835}
\newcommand{\boxd}{61\,110}
\newcommand{\boxe}{209\,431}
\newcommand{\candthree}{2\,992}
\newcommand{\candtwo}{636}
\newcommand{\dimwin}{2.4123}
\newcommand{\fdbubfour}{192}
\newcommand{\fdbubsix}{1\,604}
\newcommand{\fdbubthree}{704}
\newcommand{\fdfullthree}{5.2\cdot10^{13}}
\newcommand{\fdfulltwo}{7.5\cdot10^{9}}
\newcommand{\fdsecactive}{10\,026\,973}
\newcommand{\fdsecactiveone}{24\,892}
\newcommand{\fdsecbub}{2\,501}
\newcommand{\fdsecfrozen}{846\,778}

\newcommand{\fdsecpieces}{29\,617}
\newcommand{\fdsectiles}{10\,876\,252}
\newcommand{\fdsectime}{64}
\newcommand{\fdvactive}{37.7}
\newcommand{\fdvactivestart}{49.5}
\newcommand{\fdvbub}{3.5}
\newcommand{\fdvcore}{0.31}
\newcommand{\fdvfour}{0.19}
\newcommand{\fdvfrozen}{58.8}
\newcommand{\fdvfrozenone}{50.2}
\newcommand{\fdvfrozenstart}{50.2}
\newcommand{\fdvfrozentwo}{8.6}
\newcommand{\fdvsix}{0.09}
\newcommand{\fdvthree}{2.9}
\newcommand{\itercounts}{165,  2729,  43933,  708417}
\newcommand{\itertimes}{1.5}
\newcommand{\mindthree}{0.0015}
\newcommand{\mindtwo}{0.0017}
\newcommand{\nbubthree}{392}
\newcommand{\nbubtwo}{212}

\newcommand{\npaths}{20\,000}
\newcommand{\nsamples}{520\,000}
\newcommand{\ofcbub}{71}

\newcommand{\ofcfrac}{0.918}
\newcommand{\ofctiles}{1\,100}
\newcommand{\ofdbub}{43\,933}

\newcommand{\ofdfrac}{0.874}
\newcommand{\ofdtiles}{708\,417}
\newcommand{\ofgdbub}{13\,207}

\newcommand{\ofgdfrac}{0.868}
\newcommand{\ofgdtiles}{213\,409}
\newcommand{\ofgsdbub}{17\,587}

\newcommand{\ofgsdfrac}{0.827}
\newcommand{\ofgsdtiles}{283\,809}
\newcommand{\ofhbub}{183}

\newcommand{\ofhfrac}{0.868}
\newcommand{\ofhtiles}{2\,934}
\newcommand{\ofsbub}{219}

\newcommand{\ofsfrac}{0.815}
\newcommand{\ofstiles}{3\,520}
\newcommand{\ofssdbub}{2\,197}

\newcommand{\ofssdfrac}{0.927}
\newcommand{\ofssdtiles}{35\,209}
\newcommand{\ofwbub}{367}

\newcommand{\ofwfrac}{0.771}
\newcommand{\ofwtiles}{5\,940}
\newcommand{\pada}{1.5}

\newcommand{\padh}{3.5}

\newcommand{\perrf}{38.0}

\newcommand{\pitf}{11\,427\,701}

\newcommand{\ptgf}{29.07}

\newcommand{\ptzf}{26.10}

\newcommand{\puge}{2.5}

\newcommand{\puze}{2.1}

\newcommand{\sieveba}{31}
\newcommand{\sievebb}{32}
\newcommand{\sievebc}{212}
\newcommand{\sievebd}{392}
\newcommand{\sievebe}{1\,832}
\newcommand{\sievebf}{4\,532}
\newcommand{\sieveca}{63}
\newcommand{\sievecb}{137}
\newcommand{\sievecc}{679}
\newcommand{\sievecd}{3\,055}
\newcommand{\sievece}{13\,673}
\newcommand{\sievecf}{59\,479}
\newcommand{\sievet}{943.6}
\newcommand{\tfthreea}{5.6}
\newcommand{\tfthreeb}{12.5}
\newcommand{\tfthreec}{14.5}
\newcommand{\tfthreed}{32.6}
\newcommand{\tftwoa}{12.0}
\newcommand{\tftwob}{27.4}
\newcommand{\tftwoc}{33.0}
\newcommand{\tftwod}{33.0}
\newcommand{\tnthreea}{146}
\newcommand{\tnthreeb}{326}
\newcommand{\tnthreec}{378}
\newcommand{\tnthreed}{847}
\newcommand{\tntwoa}{51}
\newcommand{\tntwob}{116}
\newcommand{\tntwoc}{140}
\newcommand{\tntwod}{140}

\newcommand{\vrchockbeta}{3.1926}

\newcommand{\vrchockbub}{1\,163}
\newcommand{\vrchockdim}{2.412}

\newcommand{\vrchocktiles}{55\,407}

\newcommand{\vrgsbeta}{3.1926}

\newcommand{\vrgsbub}{43\,933}
\newcommand{\vrgsdim}{2.412}

\newcommand{\vrgstiles}{708\,417}

\newcommand{\vrspikebeta}{3.3001}

\newcommand{\vrspikebub}{57\,573}
\newcommand{\vrspikedim}{2.481}

\newcommand{\vrspiketiles}{1\,030\,777}

\newcommand{\vrwedgebeta}{3.1926}

\newcommand{\vrwedgebub}{2\,527}
\newcommand{\vrwedgedim}{2.412}

\newcommand{\vrwedgetiles}{29\,775}
\newcommand{\sieveratio}{2.6}
\newcommand{\crossk}{6}
\newcommand{\winova}{68}
\newcommand{\winovb}{69}
\newcommand{\winovc}{68}

\graphicspath{{./}{../}}
\newcommand{\fig}[2][\linewidth]{\includegraphics[width=#1]{#2}}
\title{The $\Zt$-foam\\[4pt]\large Goodman-Strauss's Dodecafoam as arithmetic in the icosahedral module:
  rules without geometry, isometries as integer matrices, and the fractal window of the bubbles}
\author{Miguel E.~P\'erez Ibarra\thanks{Formerly Escuela Superior de F\'isica y
    Matem\'aticas, Instituto Polit\'ecnico Nacional, Mexico City, Mexico.
    E-mail: \texttt{supersluky@gmail.com}.}
  \and Jaziel Flores Rodr\'iguez\thanks{Formerly Escuela Superior de F\'isica y
    Matem\'aticas, Instituto Polit\'ecnico Nacional, Mexico City, Mexico.
    E-mail: \texttt{jf@saptiva.com}.}}
\date{November 2026}

\begin{document}
\maketitle

\begin{abstract}
Goodman-Strauss's Dodecafoam is a substitution that fills the regular dodecahedron with
dodecahedral bubbles of sizes $\tau^{-k}$, $\tau=(1+\sqrt5)/2$, using the seven kinds of cells into
which the twelve planes of its stellation cut it.  We rewrite it as exact arithmetic in one
$\Zt$-module of rank three, the span $\M$ of the six five-fold axes, on which the icosahedral group
acts by $120$ matrices in $\GL(3,\Zt)$.  A placed inflated copy of a cell is a union of cells iff
its facet levels lie in $\{\pm1,\pm\tau^3\}$ and it lies in the dodecahedron.  Solving this finite
linear system for all scales at once gives exactly the $22$ rules known from a geometric search; no
rule exists at any scale other than $1,\tau,\tau^2,\tau^3$.  The substitution becomes integer
arithmetic on six-tuples and reproduces the bubbles of the geometric foam exactly.  Whether a point
of $\M$ is the centre of a bubble of size $\tau^{-k}$ is decided by at most $k+3$ exact sign vectors
(its address), without building the foam, and all bubbles of one size are obtained by sieving the
points of $\M$ whose Galois conjugate satisfies $|x^\ast|\le2\tau^k$, a bound we prove with computer
assistance.  The number of bubbles of size $\tau^{-k}$ is $A\beta^k+O(|\beta'|^k)$ with
$\beta=(1+\sqrt{29})/2$.  Hence the closure of the rescaled conjugates of the bubble centres has
box dimension at most $\log_\tau\beta=\dimwin<3$, and the bubbles are not selected by any
cut-and-project window with interior.  The other foams of the 1999 article are sub-foams of this
one.  The $50\,400$ substitutions allowed by the rule list give $18\,899$ bubble sets, and the
Dodecafoam is a product set of $13\,440$ of them, to the depth tested.  In four dimensions the same
linear system, for the star cut of the $120$-cell, reproduces the $718$ clusters of a geometric
search and adds the scale $\sqrt5\,\tau^2$.
\end{abstract}


\section{Introduction}\label{sec:intro}

\paragraph{The Dodecafoam.}  Extend the twelve face planes of a regular dodecahedron $\Dz$.  They
bound its final stellation, whose convex hull is a dodecahedron $\B=\tau^3\Dz$, and they cut $\B$
into $165$ convex \emph{cells} of seven shapes: the \emph{core} $\Dz$ itself and six others, which
Goodman-Strauss \cite{GS99} calls chock, spike, small and large hat, small and large wedge
(Table~\ref{tab:cells}); they are the seven elementary cells of Kramer \cite{Kramer82}.  A
\emph{rule} (or cluster) is an identity $\lambda A+t=\bigcup_i A_i$ that writes an inflated,
translated copy of a cell $A$ as a union of cells $A_i$; a \emph{placement} of a rule is one
translation $t$ realising it.  The rule $\tau^3\Dz=\B$ is the first.  Goodman-Strauss found rules
for all six non-core shapes in which one of the pieces is a core.  Starting from the \emph{seed}
$\B$ and replacing every cell that is not a core by the scaled-down pieces of its rule, again and
again, one obtains the \emph{Dodecafoam}: the cores, once produced, are never subdivided (they are
\emph{frozen}; we call them \emph{bubbles}), and in the limit $\B$ is filled by dodecahedral bubbles
of sizes $\tau^0,\tau^{-1},\tau^{-2},\dots$\ (Figure~\ref{fig:foam}).  We say that a bubble
congruent to $\tau^{-k}\Dz$ has \emph{scale} $\tau^{-k}$ and \emph{scale index} $k$, call the set of
all bubbles of scale $\tau^{-k}$ \emph{layer} $k$, and call the foam after $d$ substitution steps the
foam of \emph{depth} $d$.

\paragraph{The question.}  In a preliminary computation (unpublished; the programs are distributed
with this article) we recovered Goodman-Strauss's rules from the seven cells by a systematic search
for placed inflated copies $\lambda A+t$ that are unions of cells, and iterated them into the foam.
That computation was exact, but it was geometry with exact coordinates: polytopes, half-spaces,
rotation matrices with entries in $\tfrac12\Zt$, a bisection on clipped volumes checked afterwards
in $\Zt$.  Here we ask whether the Dodecafoam can be written down by $\Zt$ algebra alone, without
isometries acting on polytopes and without iterating the substitution.  The answer has three parts,
and the article is organised around them.

\begin{enumerate}
\item \emph{Isometries are not needed as geometric objects.}  All vertices of the cells, all
  translations of the rules and all centres of the bubbles lie in one module $\M\subset\Zt^3$ of
  rank three over $\Zt$ (Section~\ref{sec:module}).  In the basis of $\M$ the icosahedral group is
  a group of $120$ matrices in $\GL(3,\Zt)$, and a group element is determined by the permutation
  it induces on the twelve face normals.  Cells are sign vectors over the twelve cutting planes,
  and the group permutes the signs (Section~\ref{sec:cells}).
\item \emph{The rules need no polytope computation.}  A placed copy $\lambda A+t$ is a union of
  cells iff each of its facet levels $\lambda c+e\cdot t$ belongs to the set $\{\pm1,\pm\tau^3\}$
  of levels of the arrangement and the copy lies in $\B$ (Theorem~\ref{thm:rules}).  Four facets
  determine $\lambda$ and $t$, and the others are checked.  This finite problem over $\Q(\sqrt5)$
  has, for all $\lambda>0$ at once, exactly the $22$ rules of the geometric search, at the scales
  $1,\tau,\tau^2,\tau^3$ and no others (Section~\ref{sec:rules}).
\item \emph{Iteration can be avoided for membership, not for a closed form.}  The substitution
  itself is integer arithmetic on six-tuples (Section~\ref{sec:subst}).  Whether a given point of
  $\M$ is the centre of a bubble of scale $\tau^{-k}$ is decided by descending its \emph{address},
  at most $k+3$ exact sign vectors, without constructing anything else
  (Section~\ref{sec:address}).  Combined with a bound on the Galois conjugates of the bubble
  centres, this gives a sieve that produces each layer of bubbles directly
  (Figure~\ref{fig:sieve}).  A closed form in the sense of cut-and-project sets, bubble centres
  $=\{x\in\M: x^\ast\in\Om\}$ for a window $\Om$ with interior, does not exist: the closure of the
  conjugates is a set of box dimension at most $\dimwin$ and of Lebesgue measure zero, and the
  bubbles of scale $\tau^{-k}$ become sparse as $k$ grows (Section~\ref{sec:window}).
\end{enumerate}

The quantitative statement behind the third part is the following.  The number $b_k$ of bubbles of
scale $\tau^{-k}$ in the foam grown from the seed satisfies $b_{k+3}=2b_{k+2}+6b_{k+1}-7b_k$, so it
grows like $\beta^k$ with $\beta=(1+\sqrt{29})/2$, a root of $x^2-x-7$ that does not lie in
$\Q(\sqrt5)$ (Theorem~\ref{thm:beta}).  Since $\beta<\tau^3$, the bubbles of scale $\tau^{-k}$ fill
a fraction of $\B$ that tends to zero, and $\log_\tau\beta=\dimwin$ bounds the dimension of the
window (Theorem~\ref{thm:dim}).

Section~\ref{sec:compare} compares the arithmetic construction with the construction of 1999 and
measures both, Section~\ref{sec:other} applies the same arithmetic to the other foams of that
article and to the other substitutions that the rule list allows, and Section~\ref{sec:4d} solves
the same linear system in four dimensions.

\paragraph{Related work.}  The modules of Section~\ref{sec:module} are standard objects of the
theory of icosahedral quasicrystals: the $\Z$-span of the six five-fold axes is the primitive
icosahedral module, with inflation factor $\tau^3$, and its $\Zt$-span is the body-centred one
\cite{BaakeGrimm12,Huck09}.  The icosahedral tilings obtained from them by projection, such as the
Ammann--Kramer--Neri tiling, are model sets with a polyhedral window \cite{BaakeGrimm12,Moody00}.
The Dodecafoam is a different kind of object, a substitution with infinitely many tile sizes, and
the point of Section~\ref{sec:window} is that its bubbles are not a model set of that kind.  The
set $\Om$ of that section is the limit set of a graph-directed system of similarities in the sense
of Mauldin and Williams \cite{MauldinWilliams88}; we only use the elementary covering bound for its
dimension, which we prove directly.

\paragraph{Computations.}  The results of Sections~\ref{sec:module}--\ref{sec:window} are computed
by \texttt{tools/zfoam.py} and the small programs listed in Section~\ref{sec:discussion}, in exact
arithmetic over $\Zt$ (pairs of rationals $a+b\tau$).  Floating point is used only in the
Monte-Carlo sampling of the window, in the evaluation of the norms in
Proposition~\ref{prop:cone}, in the enumeration of the candidates of the sieve (with a tolerance
that can only enlarge the set), and in the renderer.  Figures~\ref{fig:foam}, \ref{fig:sieve} and
\ref{fig:other} are rendered from the integer tiles with the Blender pipeline of \cite{SDR}; the
four-dimensional figures of Section~\ref{sec:4d} are computed with the geometric substitution
program.

\section{The module, the group, the units}\label{sec:module}

\paragraph{Coordinates.}  The dodecahedron $\B$ has vertices $(\pm\tau,\pm\tau,\pm\tau)$ and the
cyclic permutations of $(0,\pm\tau^2,\pm1)$.  Its twelve outer face normals are the cyclic
permutations of $(\pm1,\pm\tau,0)$; we number them $e_1,\dots,e_{12}$ in the order of the catalogue
of \cite{SDR}.  They are the vertices of an icosahedron, six antipodal pairs $\pm e$ (the six
five-fold axes), with $|e_i|^2=2+\tau$.  The faces of $\B$ lie in the planes $e_i\cdot x=\tau^3$
($\tau^3=1+2\tau$) and the twelve cutting planes are $e_i\cdot x=1$, $i=1,\dots,12$.  The
central cell is the core $\Dz=\{x:e_i\cdot x\le1,\ i=1,\dots,12\}$, and $\B=\tau^3\Dz$.

\begin{definition}
The \emph{icosahedral module} is $\M=\Zt e_1+\dots+\Zt e_{12}\subset\Zt^3$.  We write
$\Pm=\Z e_1+\dots+\Z e_{12}$ for the $\Z$-span of the same vectors.
\end{definition}

\begin{proposition}\label{prop:module}
\begin{enumerate}
\item $\M$ is a free $\Zt$-module of rank $3$ with basis $g_1=(1,\tau,0)$, $g_2=(0,1,\tau)$,
  $g_3=(\tau,0,1)$; $\det(g_1,g_2,g_3)=2\tau^2$ and $[\Zt^3:\M]=4$.
\item $\Pm$ has index $2$ in $\M$.  It is invariant under multiplication by $\tau^3$ but not by
  $\tau$ or $\tau^2$.  $\M=\Pm+\tau\Pm$ is the smallest $\tau$-invariant module containing the
  five-fold axes, and $\M=\Pm\cup\bigl(\Pm+(1,1,1)\bigr)$.
\item The vertices of all $165$ cells lie in $\M$.
\end{enumerate}
\end{proposition}

\begin{proof}
(1) The three vectors are $\Zt$-independent since their determinant $1+\tau^3=2\tau^2$ is not zero,
and each $e_i$ is a $\Zt$-combination of them (for example $(0,1,-\tau)=g_1+(1-\tau)g_2+
(1-\tau)g_3$), so they span $\M$.  The index is the absolute value of the norm of the determinant,
$|N(2\tau^2)|=4$.  (2) The $6\times6$ integer matrix of six pairwise non-antipodal normals in the
$\Z$-basis $g_j,\tau g_j$ of $\M$ has determinant $\pm2$.  Solving the integer linear systems
exactly gives $\tau e_1,\tau^2e_1\notin\Pm$ and $\tau^3e_1\in\Pm$; the icosahedral group preserves
$\Pm$, commutes with multiplication by $\tau^3$ and is transitive on the normals, so
$\tau^3\Pm\subset\Pm$.  $\M=\Pm+\tau\Pm$ because $\Zt=\Z+\Z\tau$.  Finally
$(1,1,1)=\tau^{-2}(g_1+g_2+g_3)\in\M$, and $(1,1,1)\notin\Pm$: a combination $\sum n_ie_i$ of
$(1,\pm\tau,0)$, $(\pm\tau,0,1)$, $(0,1,\pm\tau)$ with first coordinate $1$ and second coordinate
free of $\tau$ would need $n_1+n_2=1$ and $n_1-n_2=0$.  (3) The vertices of the seven
representative cells of the exact catalogue \cite{SDR} are in $\M$ by direct check, every cell is
the image of a representative under the group, and the group preserves $\M$ (next paragraph).
\end{proof}

Part (2) is the standard fact that the primitive icosahedral module has the inflation factor
$\tau^3$, while the body-centred module $\Pm\cup(\Pm+(1,1,1))$ is a $\Zt$-module
\cite[Section~4.5]{BaakeGrimm12}, \cite[Section~3]{Huck09}.  It matches the two inflation factors
that occur here: Kramer's seven cells inflate by $\tau^3$ \cite{Kramer82}, and Goodman-Strauss's
rules use the scales $\tau$ and $\tau^2$ as well, with translations that lie in $\M$.

\paragraph{The group.}  An isometry of $\B$ permutes the twelve normals.  Choosing three independent
normals and prescribing their images among the $12^3$ triples with the same Gram matrix, the
program finds $120$ linear maps that permute all twelve normals: the full icosahedral group
$I_h$, $60$ rotations and $60$ improper elements.  In the standard coordinates their matrices have
entries in $\tfrac12\Zt$; since they map $\M$ onto itself, in the basis $g_1,g_2,g_3$ they lie in
$\GL(3,\Zt)$.  For example the rotation
\[
R=\tfrac12\begin{pmatrix}1-\tau&\tau&-1\\ \tau&1&\tau-1\\ 1&1-\tau&-\tau\end{pmatrix}
\quad\text{becomes}\quad
A=\begin{pmatrix}1&1&\tau-1\\ 0&1-\tau&\tau-1\\ 0&1-\tau&-1\end{pmatrix}\in\GL(3,\Zt).
\]
A group element $g$ is stored as its permutation $\pi_g$ of the twelve normals
($g\,e_i=e_{\pi_g(i)}$), which is what acts on sign vectors, and as its matrix $A_g$, which is what
acts on points of $\M$; both are integer data.  Writing $x\in\M$ as $x=\sum_j(a_j+b_j\tau)g_j$ with
six integers, $A_g$ acts as a $6\times6$ integer matrix and multiplication by $\tau^{\pm1}$ as
$(a,b)\mapsto(b,a+b)$, $(a,b)\mapsto(b-a,a)$.

\paragraph{Units and levels.}  The arrangement has two levels, $1$ for the cuts and $\tau^3$
for the faces, and their ratio is a unit of $\Zt$.  Consequently multiplication by $\tau^k$, which
maps $\M$ onto itself, maps the set of cut planes $\{e_i\cdot x=1\}$ onto $\{e_i\cdot x=\tau^k\}$,
which for $k=3$ is the set of face planes.  This is the arithmetic content of the rule
$\tau^3\Dz=\B$, and Section~\ref{sec:rules} shows that all other rules are solutions of equations
between the same two levels.

\paragraph{Galois conjugation.}  The nontrivial automorphism of $\Q(\sqrt5)$, $\tau\mapsto
\tau'=1-\tau=-1/\tau$, applied coordinate-wise, is the \emph{star map} $x\mapsto x^\ast$ on $\M$.
It is additive, $(\tau x)^\ast=-\tau^{-1}x^\ast$, and $\{(x,x^\ast):x\in\M\}$ is a lattice of rank
$6$ in $\R^3\times\R^3$ (the lift of $\M$).  If $R_g$ is the matrix of $g$ in the standard
coordinates then $(gx)^\ast=R_g^\ast x^\ast$, where $R_g^\ast$ is the conjugate matrix; it is
orthogonal, because the relation $R_g^{\mathsf T}R_g=1$ is preserved by the automorphism.  This is
the second three-dimensional representation of $I_h$.  The second factor $\R^3$ is called internal
space; the window of Section~\ref{sec:window} lives there.

\section{Cells as sign vectors}\label{sec:cells}

A point $x\in\B$ off the cuts has a sign vector $\sigma(x)\in\{-,+\}^{12}$, $\sigma_i=\sgn(e_i\cdot
x-1)$.  A \emph{cell} is the closure of the set of points with a given sign vector; two cells
overlap in at most a face.  Table~\ref{tab:cells} lists the seven classes by the sign vectors of
their representatives, with Goodman-Strauss's names and the exact volumes (sums of cones over the
exact vertices; with the multiplicities they add up to $\vol\B=16+28\tau$).  The group acts by
$\sigma_{\pi_g(i)}(gx)=\sigma_i(x)$, and the $165$ cells are the orbits of the seven
representatives; the program reproduces the orbit sizes $30,1,20,12,60,12,30$.  A cell
$g(\mathrm{rep}_c)$ is stored as the pair (class $c$, group element $g$); no translation is
involved, since the arrangement is centred at the origin.

\begin{table}[ht]\centering\small
\caption{The seven classes of cells as sign vectors over the twelve cuts $e_i\cdot x=1$ (order of
  the normals as in the catalogue of \cite{SDR}: the cyclic permutations of $(\pm1,\pm\tau,0)$).
  Names after Goodman-Strauss \cite{GS99}, who writes $D_0$ for the core (``small dodecahedron'');
  $H_1=\tau H_0$ and $W_1=\tau W_0$ up to translation.}\label{tab:cells}
\begin{tabular}{llrlr}
\toprule
class & name & orbit & sign vector & volume \\
\midrule
$C_0$ & chock & 30 & \verb|------++--++| & $(10-4\tau)/3$ \\
$\Dz$ & core & 1 & \verb|------------| & $72-44\tau$ \\
$S_0$ & spike & 20 & \verb|-------+--++| & $-6+4\tau$ \\
$H_1$ & large hat & 12 & \verb|--++--++--++| & $(14-8\tau)/3$ \\
$W_1$ & large wedge & 60 & \verb|---+--++--++| & $(-6+4\tau)/3$ \\
$H_0$ & small hat & 12 & \verb|-----------+| & $(-58+36\tau)/3$ \\
$W_0$ & small wedge & 30 & \verb|----------++| & $(26-16\tau)/3$ \\
\bottomrule
\end{tabular}
\end{table}

A \emph{facet} of a class is a triple $(e_i,c,s)$: the facet lies in the plane $e_i\cdot x=c$ with
$c\in\{1,\tau^3\}$, and $s=+1$ if the cell lies on the side $e_i\cdot x\le c$, $s=-1$ otherwise.
The facets are read off the sign vector and the face incidences of the catalogue.  Apart from the
vertices and volumes used for exact checks, this is all the geometry the rest of the article uses.

\section{The rules as a linear system}\label{sec:rules}

\begin{theorem}\label{thm:rules}
Let $A$ be a cell with facets $(e_{i_j},c_j,s_j)$, let $\lambda>0$ and $t\in\R^3$.  Then
$\lambda A+t$ is a union of cells of the arrangement iff
\begin{enumerate}
\item every level $L_j=\lambda c_j+e_{i_j}\cdot t$ lies in $\{\pm1,\pm\tau^3\}$,
\item a facet landing on a face plane of $\B$ ($L_j=\pm\tau^3$) has $\lambda A+t$ on the inner side,
\item $\lambda A+t\subset\B$.
\end{enumerate}
In that case the cells of $\lambda A+t$ are exactly those whose sign vectors satisfy one sign
constraint for each facet with $L_j=\pm1$: the sign of the cut $e_{i_j}\cdot x=1$ if $L_j=1$, of
the antipodal cut $-e_{i_j}\cdot x=1$ if $L_j=-1$, as prescribed by the side $s_j$.  If moreover
$\lambda\in\Zt$, then $t\in\M$.
\end{theorem}

Condition (2) follows from (3); it is listed because it is the test the program applies first.
Condition (3) is not implied by (1) and (2): the cut planes extend beyond $\B$, and a polytope
bounded by cuts may leave $\B$ (this happens, for instance, for placements of $\tau^2C_0$,
$\tau^3H_1$ and $\tau^2W_1$).

\begin{proof}
Put $P=\lambda A+t$, a convex polytope with nonempty interior whose facet $j$ lies in the plane
$e_{i_j}\cdot x=L_j$, on the side given by $s_j$.  If $P$ is a union of cells, then $P\subset\B$
and the boundary of $P$ is contained in the union of the $24$ planes of the cuts and faces, so
every facet plane of $P$ is one of them: this is (1) and (3), and (2) follows.  Conversely, assume
(1)--(3).  $P$ is the intersection of the half-spaces of its facets.  By (2) the half-spaces of the
facets with $L_j=\pm\tau^3$ contain $\B$, so by (3) $P$ is the intersection of $\B$ with the
half-spaces of the facets with $L_j=\pm1$, each bounded by a cut.  A cell lies on one side of every
cut, hence $P$ is the union of the cells on the prescribed side of each of these cuts, which is the
description by sign constraints.  Finally let $\lambda\in\Zt$ and let $v$ be a vertex of $A$.  Then
$\lambda v\in\M$ by Proposition~\ref{prop:module}(3), and $\lambda v+t$ is a vertex of $P$, hence
an extreme point of a cell contained in $P$, that is, a vertex of a cell, which lies in $\M$.  So
$t\in\M$.
\end{proof}

The theorem turns the geometric cluster search into a finite computation over $\Q(\sqrt5)$.  For a
given scale $\lambda$, three facets with independent normals and the $4^3$ choices of their target
levels determine $t$ by Cramer's rule.  The remaining facets, the orientation and the inclusion in
$\B$ are checked exactly, and the composition is read off the $165$ sign vectors and verified
against the volume identity $\sum_c n_c\vol(c)=\lambda^3\vol(A)$.  The scale need not be given:
every cell has four facets whose vectors $(c_j,e_{i_j})\in\Q(\sqrt5)^4$ are independent, and the
$4^4$ choices of their target levels determine the pair $(\lambda,t)$.  Both computations take a
few seconds (\texttt{zfoam.py rules} for the scales $\tau^k$, $0\le k\le6$, $\sqrt5\,\tau^k$,
$-2\le k\le2$, and $2$; \texttt{allscales.py} for all $\lambda>0$).

\begin{result}[the complete rule set]\label{res:rules}
For $\lambda>0$ and $t\in\R^3$ arbitrary, the solutions of Theorem~\ref{thm:rules} have
$\lambda\in\{\tau^{-1},1,\tau,\tau^2,\tau^3\}$ and $t\in\M$.  Those with $\lambda\le1$ and a single
cell are the identities, one translation that carries the representative of $W_1$ onto another cell
of its class, and the inverses $\tau^{-1}H_1=H_0$, $\tau^{-1}W_1=W_0$ of two rules below.  All the
others are the $22$ rules of Table~\ref{tab:rules}, which are exactly the rules found by the
geometric search, with the same numbers of placements.  Goodman-Strauss's Dodecafoam rules
\cite{GS99} are among them: $C_0=\Dz+4H_0+W_0$, $\tau S_0=\Dz+S_0+6H_0+6W_0$,
$\tau H_1=\tau^2H_0=\Dz+6H_0+5W_0$, $\tau^2W_1=\tau^3W_0=\Dz+2S_0+8H_0+11W_0$, and
$\tau^3\Dz=\B$.
\end{result}

\begin{table}[ht]\centering\small
\caption{The $22$ rules as solutions of the linear system, with one translation each (coordinates
  in $\Zt^3$) and the number of placements, i.e.\ of translations $t$ for which the inflated
  representative of the class has this composition.}
\label{tab:rules}
\begin{tabular}{llrl}
\toprule
scale & rule & placements & $t$ \\
\midrule
$1$ & $C_0=\Dz+4H_0+W_0$ & 1 & $(0,0,2)$ \\
$\tau$ & $\tau C_0=3C_0+2S_0+3W_1+W_0$ & 2 & $(0,\tau,\tau^2)$ \\
$\tau^3$ & $\tau^3\Dz=30C_0+\Dz+20S_0+12H_1+60W_1+12H_0+30W_0$ & 1 & $0$ \\
$\tau$ & $\tau S_0=C_0+H_1+3W_1$ & 3 & $(0,\tau,\tau^2)$ \\
$\tau$ & $\tau S_0=\Dz+S_0+6H_0+6W_0$ & 1 & $(\tau^{-1},0,\tau)$ \\
$\tau^2$ & $\tau^2S_0=3C_0+\Dz+5S_0+3W_1+9H_0+15W_0$ & 1 & $(\tau,0,\tau^3)$ \\
$\tau$ & $\tau H_1=C_0+2W_1$ & 5 & $(0,2+\tau,\tau^2)$ \\
$\tau$ & $\tau H_1=\Dz+6H_0+5W_0$ & 1 & $(0,\tau,\tau^2)$ \\
$1$ & $W_1=H_0+2W_0$ & 2 & $(0,0,2)$ \\
$\tau$ & $\tau W_1=S_0+H_0+3W_0$ & 2 & $(0,\tau,\tau^2)$ \\
$\tau$ & $\tau W_1=H_1+2W_1$ & 2 & $(\tau,\tau^2,2\tau)$ \\
$\tau^2$ & $\tau^2W_1=C_0+2H_1+6W_1$ & 1 & $(\tau,0,\tau^3)$ \\
$\tau^2$ & $\tau^2W_1=\Dz+2S_0+8H_0+11W_0$ & 1 & $(1,\tau,2\tau^2)$ \\
$\tau$ & $\tau H_0=H_1$ & 1 & $(0,\tau,\tau^2)$ \\
$\tau^2$ & $\tau^2H_0=C_0+2W_1$ & 5 & $(\tau,0,1)$ \\
$\tau^2$ & $\tau^2H_0=\Dz+6H_0+5W_0$ & 1 & $(0,1,\tau)$ \\
$\tau$ & $\tau W_0=H_0+2W_0$ & 2 & $(\tau^{-1},0,\tau)$ \\
$\tau$ & $\tau W_0=W_1$ & 2 & $(0,\tau,\tau^2)$ \\
$\tau^2$ & $\tau^2W_0=H_1+2W_1$ & 2 & $(\tau,0,1)$ \\
$\tau^2$ & $\tau^2W_0=S_0+H_0+3W_0$ & 2 & $(0,-1,\tau)$ \\
$\tau^3$ & $\tau^3W_0=\Dz+2S_0+8H_0+11W_0$ & 1 & $(0,0,2\tau)$ \\
$\tau^3$ & $\tau^3W_0=C_0+2H_1+6W_1$ & 1 & $(0,0,2\tau^2)$ \\
\bottomrule
\end{tabular}
\end{table}

\paragraph{The arithmetic behind the rules.}  Condition (1) says that for each facet the shift
$e_i\cdot t$ equals $L-\lambda c$ with $L\in\{\pm1,\pm\tau^3\}$ and $c\in\{1,\tau^3\}$; replacing
$e_i$ by its antipode changes all signs.  Up to sign, for $\lambda=\tau$ the admissible shifts are
$\tau^{-1}$, $\tau^2$, $1+3\tau$, $3+3\tau$, $3+5\tau$; for $\lambda=\tau^2$ they are $\tau$,
$2+\tau$, $2+3\tau$, $2+5\tau$, $4+5\tau$, $4+7\tau$; for $\lambda\in\{1,\tau^3\}$ they are even:
$0$, $2$, $2\tau$, $2+2\tau$, \dots\ The rules are the simultaneous solutions of these conditions on
the facets of one cell, and the identities $\tau^3-\tau=\tau^2$, $\tau^3+1=2\tau^2$,
$\tau^3-1=2\tau$ between units are what make them solvable.  Conditions (1) and (2) alone still
have solutions beyond $\tau^3$, for instance $\tau^4W_0+(1,-1,1+2\tau)$; what excludes the larger
scales is the inclusion in $\B$.  For $\tau^k$ with $k\ge6$ this is visible on the volumes alone:
the smallest cell has $\tau^{18}\vol(W_0)=(178+288\tau)/3>16+28\tau=\vol\B$.

\section{The substitution as integer arithmetic}\label{sec:subst}

A tile of the foam is a quadruple $(c,k,g,u)$: class $c$, scale index $k$, group element $g$,
position $u\in\M$, representing the polytope $\tau^{-k}g(\mathrm{rep}_c)+u$.  One rule per class is
chosen as Goodman-Strauss does \cite{GS99}; in our terms it is, for each class, the rule that
contains a core at the smallest scale, and the core itself is expanded only as the seed.  The
inflation exponents $m_c$ are $m=0$ for $C_0$, $1$ for $S_0,H_1$, $2$ for $W_1,H_0$, $3$ for $W_0$,
and $3$ for the seed.  If the chosen rule is $\tau^{m_c}\mathrm{rep}_c+t_c=
\bigcup_i g_i(\mathrm{rep}_{c_i})$ then
\[
\tau^{-k}g(\mathrm{rep}_c)+u=\bigcup_i\;\tau^{-(k+m_c)}\,(g\,g_i)(\mathrm{rep}_{c_i})
+\bigl(u-\tau^{-(k+m_c)}g(t_c)\bigr),
\]
so the children of $(c,k,g,u)$ are the tiles
\begin{equation}\label{eq:children}
(c_i,\;k+m_c,\;g\,g_i,\;u-\tau^{-(k+m_c)}A_g t_c).
\end{equation}
The composition $g\,g_i$ is a product of permutations, $A_gt_c$ a $6\times6$ integer matrix applied
to six integers, and $\tau^{-(k+m_c)}$ the integer map $(a,b)\mapsto(b-a,a)$ iterated.  Cores are
frozen: they are the bubbles.  The seed is $(\Dz,-3,\mathrm{id},0)$, that is, $\B=\tau^3\Dz$.

\begin{lemma}\label{lem:inM}
Every tile $(c,k,g,u)$ produced from the seed has $u\in\M$, and a bubble $(\Dz,k,g,u)$ is the
dodecahedron $\tau^{-k}\Dz+u$ with centre $u$.  All children of a tile have the same position.
\end{lemma}

\begin{proof}
The seed has $u=0$.  In \eqref{eq:children}, $t_c\in\M$ by Result~\ref{res:rules}, $A_g$ maps $\M$
to itself, and so does $\tau^{-1}$ because $\tau$ is a unit; the new position does not depend on
$i$.  Finally $g(\Dz)=\Dz$, which is centred at the origin.
\end{proof}

\begin{result}[the foam in $\M$]\label{res:iterate}
From the seed, the integer iteration produces \itercounts{} tiles at depths $1,\dots,4$ in
\itertimes{} seconds.  The sets of pairs (scale, centre) of the bubbles coincide exactly with those
of the geometric iteration at depth $3$ ($2\,729$ bubbles) and depth $4$ ($43\,933$ bubbles).  At
depth $4$ the bubbles have the scales $\tau^0,\dots,\tau^{-9}$, in numbers
$31,32,212,392,1812,4212,8892,12960,11760,3630$.
\end{result}

\begin{figure}[p]\centering
\fig[0.49\linewidth]{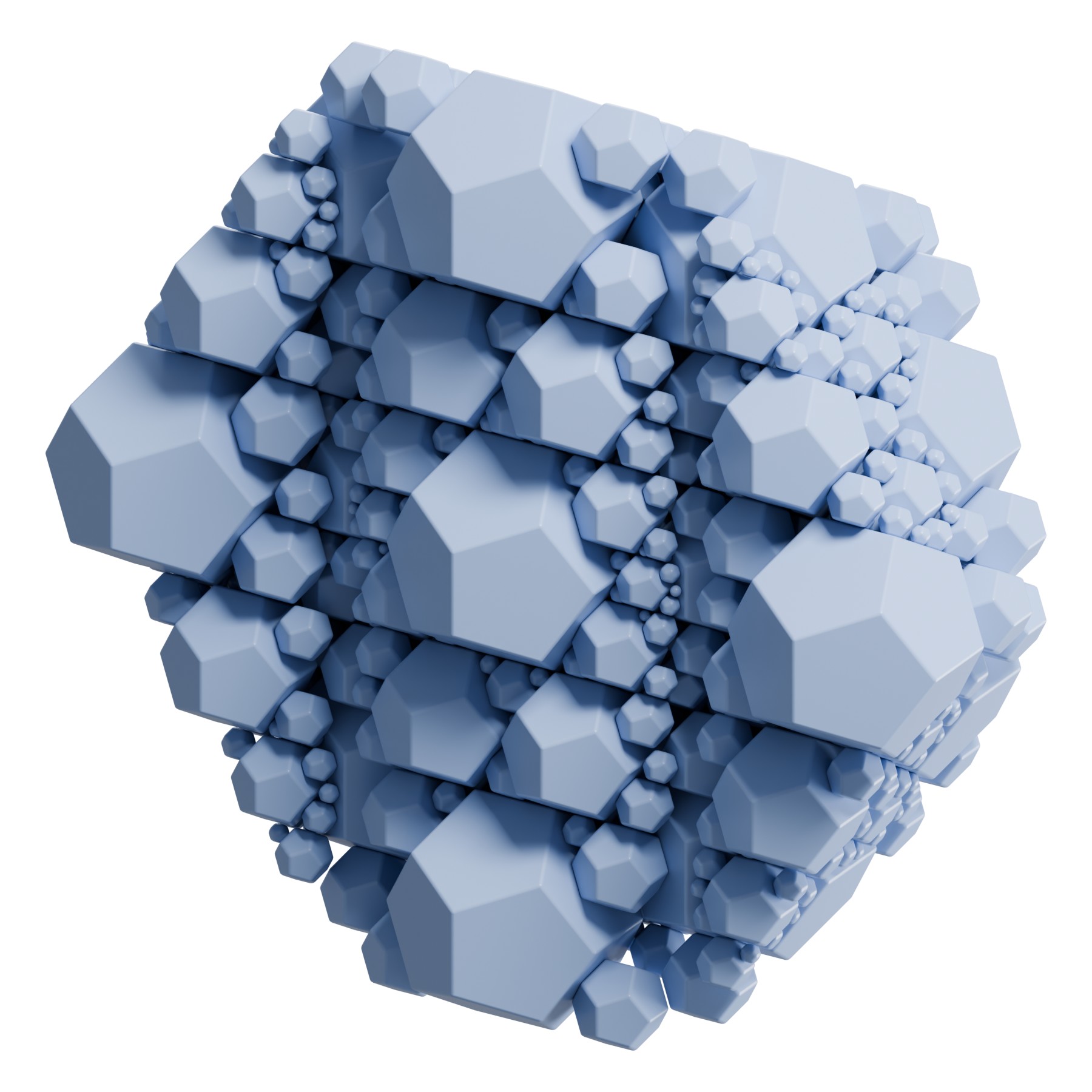}\hfill
\fig[0.49\linewidth]{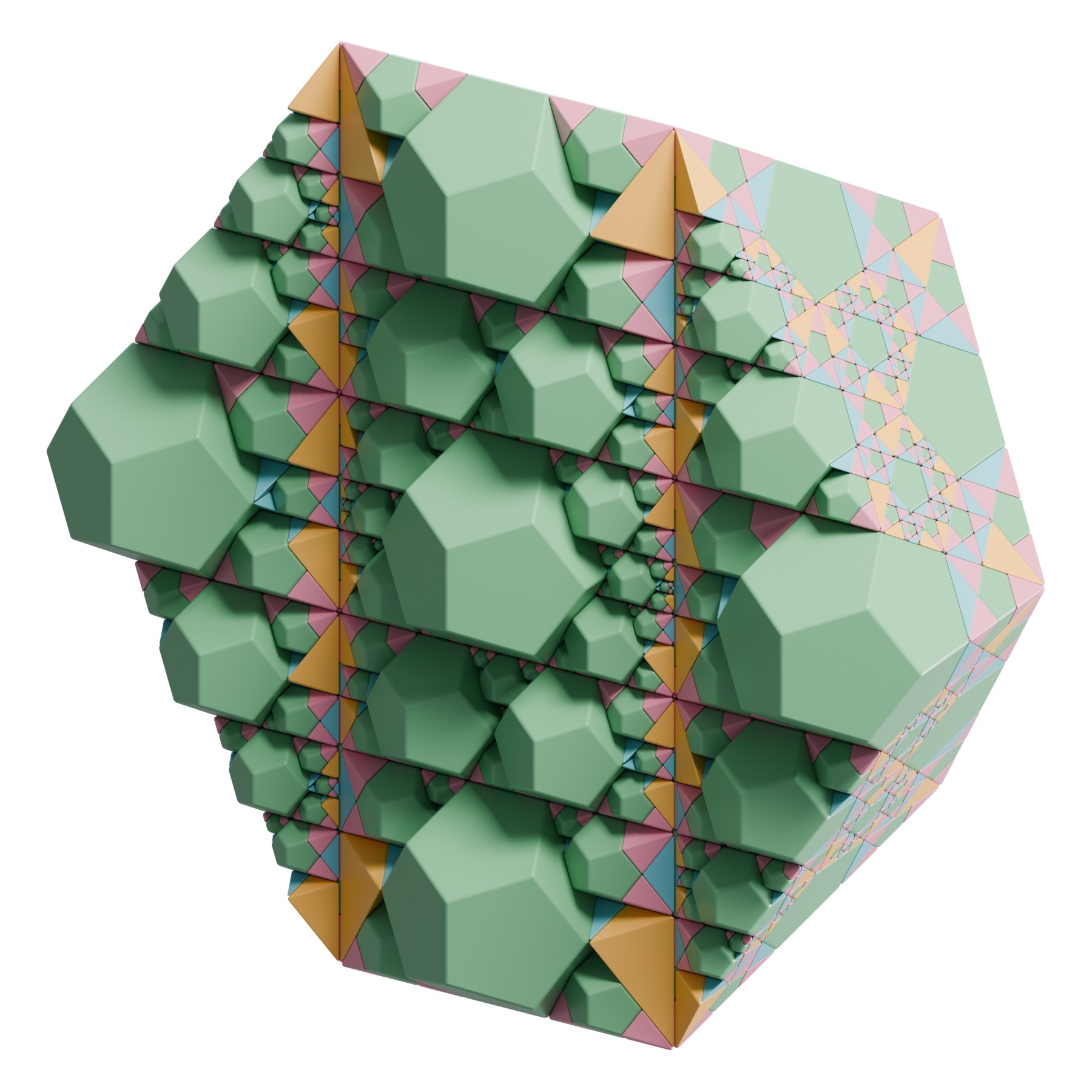}\\[4pt]
\fig[0.6\linewidth]{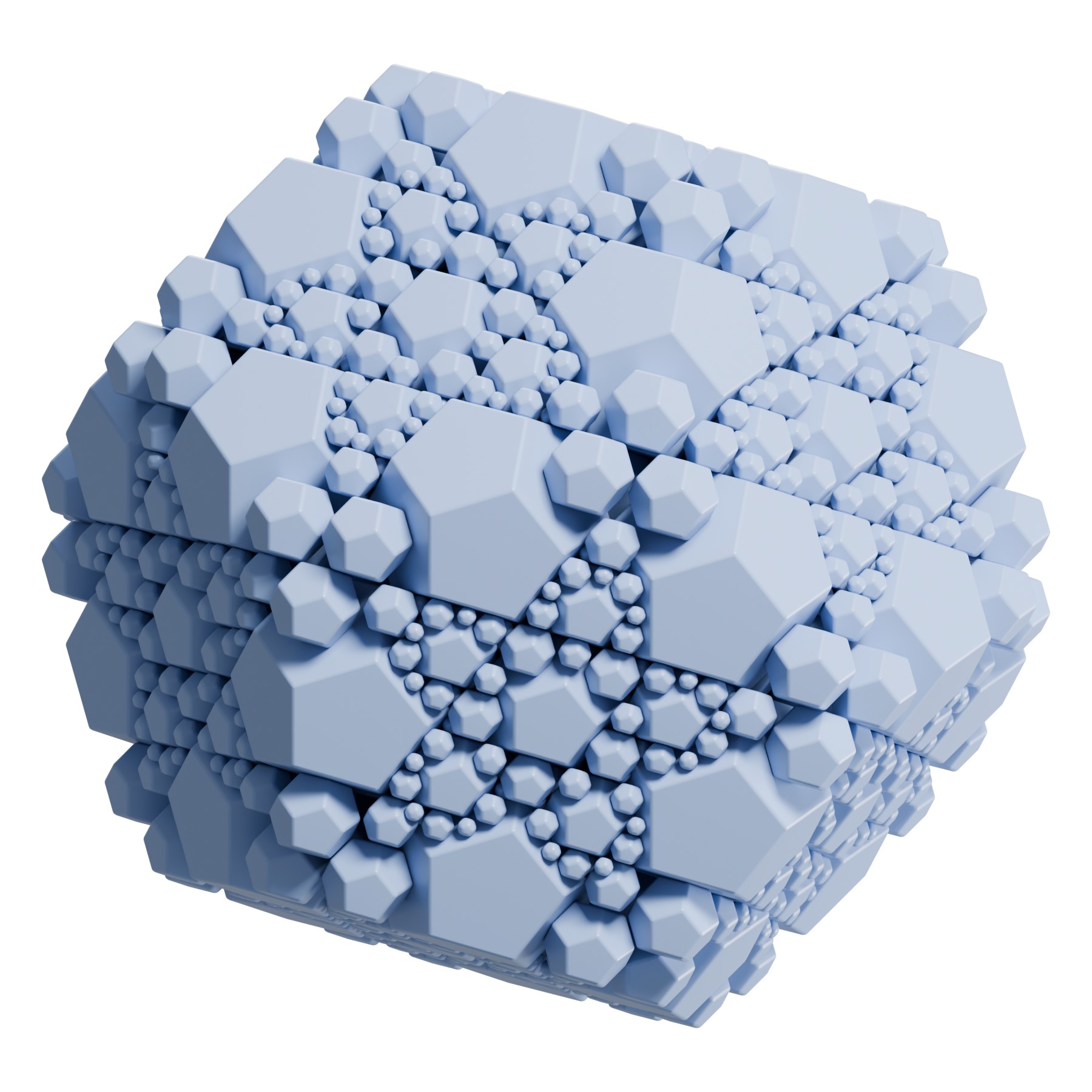}
\caption{The Dodecafoam after three substitution steps, computed in the module $\M$ (integer
  six-tuples and permutations; floating point only in the renderer).  Top left: the $2\,729$
  dodecahedral bubbles of scales $\tau^0,\dots,\tau^{-6}$, cut in half.  Top right: all
  $43\,933$ tiles, coloured by class (bubbles blue; spikes, hats and wedges in the other
  colours), cut in half.  Bottom: the bubbles, whole.  The bubbles are the same as those of the
  geometric iteration (Result~\ref{res:iterate}).}\label{fig:foam}
\end{figure}

\begin{figure}[tp]\centering
\fig[0.49\linewidth]{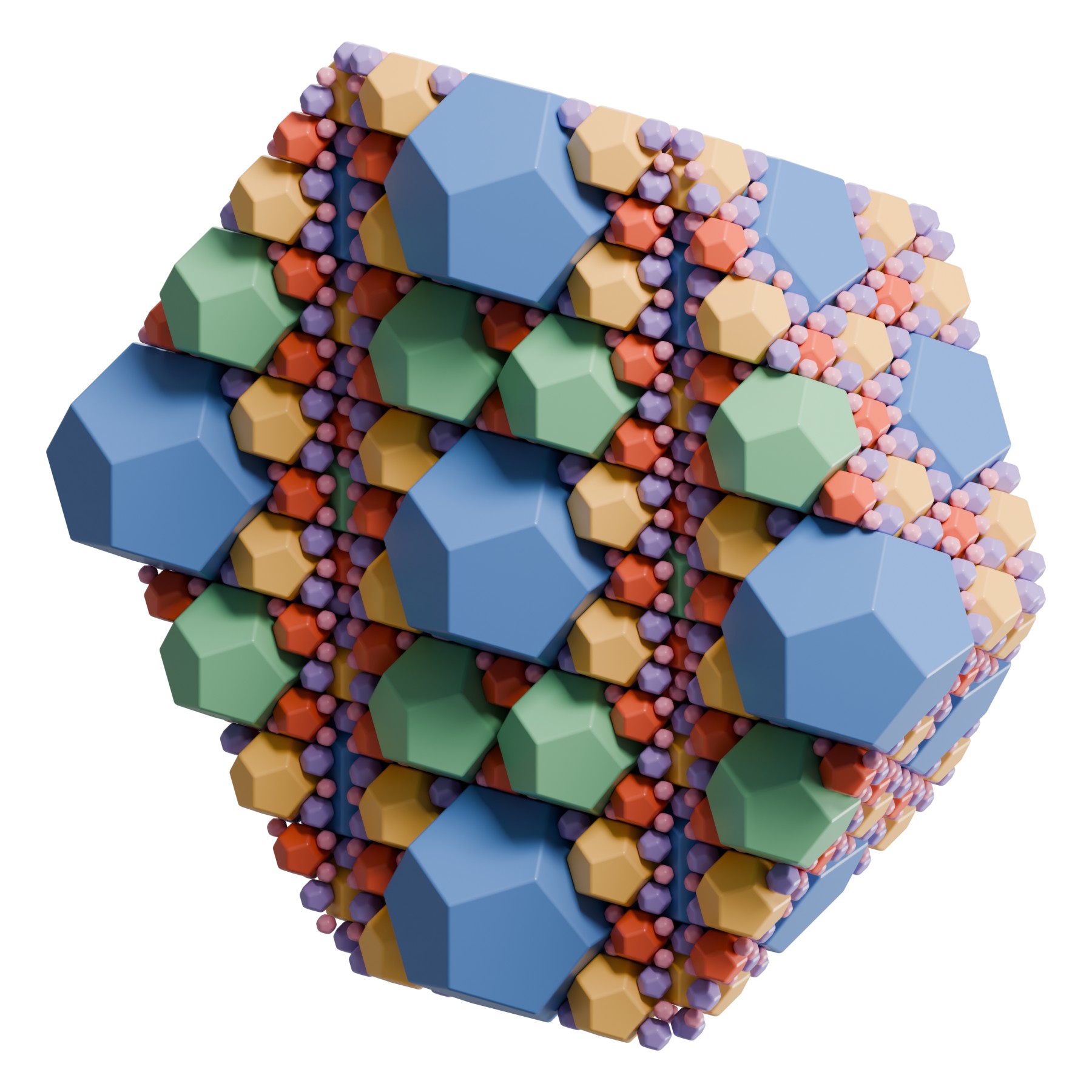}\hfill
\fig[0.49\linewidth]{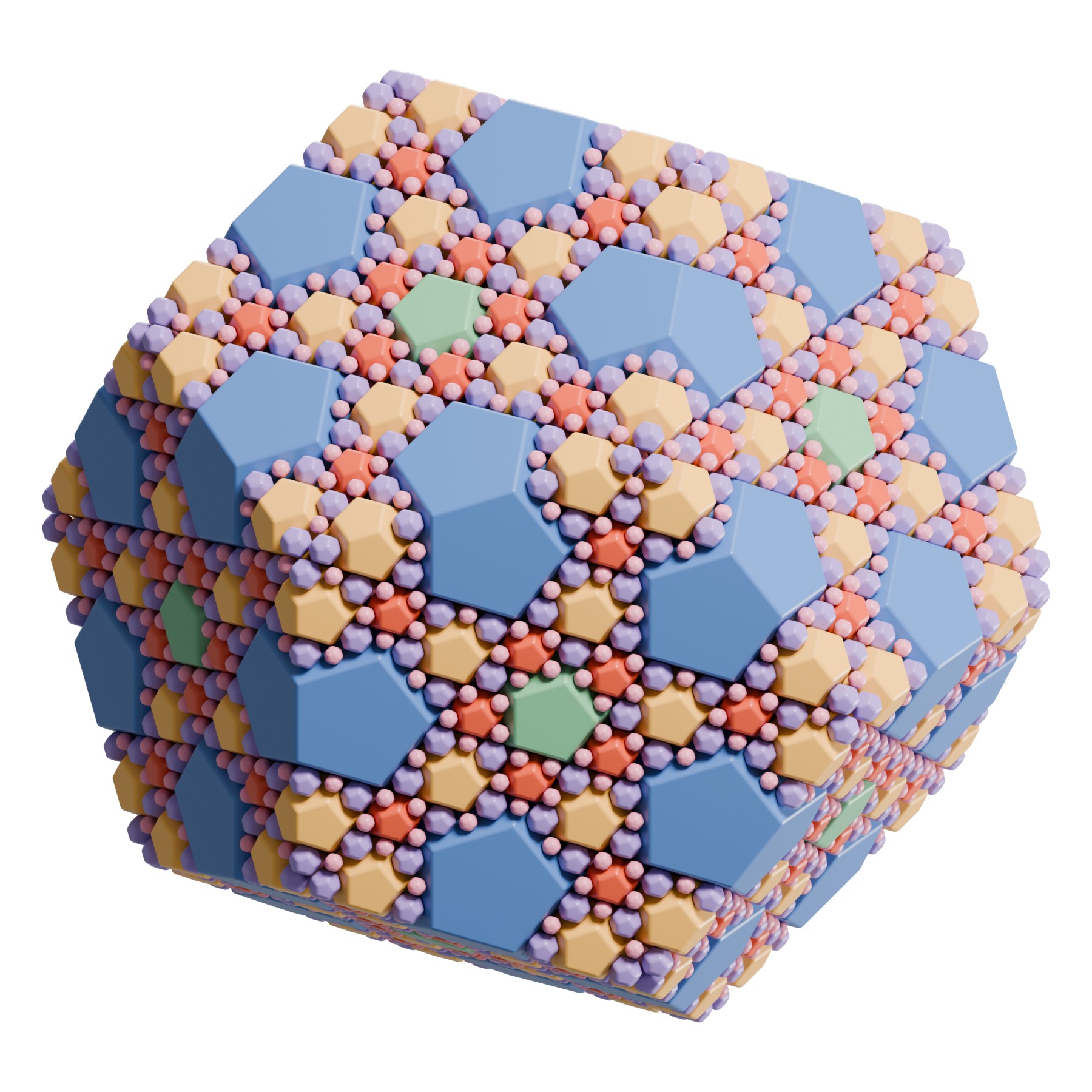}
\caption{The Dodecafoam produced without iterating it: the bubbles of scales $\tau^0,\dots,
  \tau^{-5}$ (coloured by scale), obtained as the points of the Galois cone whose address ends at
  a core of the right scale (Result~\ref{res:sieve}).  Left: cut in half.  Right: whole.  Unlike
  in the depth-$4$ foam, the layers $4$ and $5$ are complete.}
\label{fig:sieve}
\end{figure}

Each chosen rule contains exactly one core, so every tile that is expanded contributes one new
bubble: the number of bubbles at depth $d+1$ equals the number of tiles at depth $d$, which is why
$2\,729$ and $43\,933$ appear both as tile counts and as bubble counts.  Figure~\ref{fig:foam} is
the depth-$3$ foam rendered from the integer tiles: the polygons of a tile have the exact vertices
$\tau^{-k}A_gv+u$, with $v$ the vertices of the representative, and are converted to floating point
only inside the renderer.

\paragraph{Counting.}  Let $a_k(c)$ be the number of active (non-core) tiles of class $c$ and
scale index $k$ produced in the whole history of the foam grown from the seed, and $b_k$ the number
of bubbles of scale index $k$.  If $n_{cc'}$ is the number of cells of class $c'$ in the rule of
class $c$, then \eqref{eq:children} gives $a_{k}(c')=\sum_c n_{cc'}\,a_{k-m_c}(c)$ and
$b_k=\sum_c n_{c\Dz}\,a_{k-m_c}(c)$, starting from the seed $a_{-3}(\Dz)=1$.  The classes $C_0$,
$H_1$, $W_1$ occur only at $k=0$, among the cells of the seed; all later tiles are of the classes
$S_0,H_0,W_0$, with delays $1,2,3$.  Hence for $k\ge3$
\begin{gather*}
\begin{pmatrix}a_k(S_0)\\ a_k(H_0)\\ a_k(W_0)\end{pmatrix}
=\begin{pmatrix}1\\6\\6\end{pmatrix}a_{k-1}(S_0)
+\begin{pmatrix}0\\6\\5\end{pmatrix}a_{k-2}(H_0)
+\begin{pmatrix}2\\8\\11\end{pmatrix}a_{k-3}(W_0),\\
b_k=a_{k-1}(S_0)+a_{k-2}(H_0)+a_{k-3}(W_0),
\end{gather*}
with the initial values $(a_k(S_0),a_k(H_0),a_k(W_0))=(20,132,60)$, $(20,192,180)$,
$(140,1392,1440)$ for $k=0,1,2$.

\begin{theorem}\label{thm:beta}
For the rules chosen above, the numbers of bubbles of scale $\tau^{-k}$ in the foam grown from
$\tau^3\Dz$ are
\begin{gather*}
b_0,\dots,b_{6}=\bka,\ \bkb,\ \bkc,\ \bkd,\ \bke,\ \bkf,\ \bkg,\\
b_7,\dots,b_{12}=\bkh,\ \bki,\ \bkj,\ \bkk,\ \bkl,\ \bkm,
\end{gather*}
and $b_{k+3}=2b_{k+2}+6b_{k+1}-7b_k$ for all $k\ge1$.  Consequently
$b_k=A\beta^k+B\beta'^{\,k}+\tfrac{44}7$ for $k\ge1$, with $\beta=(1+\sqrt{29})/2=\betaval$,
$\beta'=(1-\sqrt{29})/2$ and $A>0$; in particular $b_{k+1}/b_k\to\beta$ and
$b_k\,\tau^{-3k}\to0$.
\end{theorem}

\begin{proof}
The values are computed from the recursion in exact integer arithmetic.  The delay system has the
characteristic polynomial
\[
\det\bigl(x^3-x^2M_1-xM_2-M_3\bigr)=x^3\,(x-1)(x^2-x-7)(x^3+x^2+2x-2),
\]
where $M_1,M_2,M_3$ are the $3\times3$ matrices whose only nonzero columns are the three columns
displayed above (the factorisation is checked by exact polynomial multiplication), so $(b_k)$
satisfies a linear recurrence of order $9$ for $k\ge3$.  The recurrence of order $3$ in the
statement, with characteristic polynomial $(x-1)(x^2-x-7)$, is verified for $1\le k\le34$; the
defect $b_{k+3}-2b_{k+2}-6b_{k+1}+7b_k$ satisfies the recurrence of order $9$ and vanishes on more
than nine consecutive terms, hence identically.  The closed form follows, the constant being
determined by $b_1,b_2,b_3$.  Since $b_k$ is positive and unbounded and $|\beta'|=2.1926<\beta$,
the coefficient $A$ is positive.  The roots of the cubic factor ($0.651$ and a complex pair of
modulus $1.753$) do not occur in $b_k$.  Finally $\beta<\tau^3=4.2361$.
\end{proof}

\begin{corollary}\label{cor:density}
The bubbles of scale $\tau^{-k}$ fill the fraction $b_k\,\tau^{-3k}\vol(\Dz)/\vol(\B)=
b_k\,\tau^{-3k-9}$ of $\B$, which tends to zero like $(\beta/\tau^3)^k$.  Equivalently, in the foam
rescaled by $\tau^k$, in which these bubbles have unit size and fill the dodecahedron $\tau^k\B$,
the number of unit bubbles per unit volume tends to zero.
\end{corollary}

That the fractions tend to zero is automatic, since the bubbles are disjoint and the fractions add
up to at most $1$; this forces $\beta\le\tau^3$ for any foam of this kind.  The content of the
corollary is the rate.

The ratios $b_{k+1}/b_k$ oscillate around $\beta$ because $\beta'<0$.  The depth-$4$ counts of
Result~\ref{res:iterate} agree with $b_k$ for $k\le3$ and fall short for $k\ge4$, since after four
steps some tiles that still contain bubbles of scale $\tau^{-4}$ have not been expanded.

A second exponent describes how fast the foam fills $\B$.  With $P(s)=sM_1+s^2M_2+s^3M_3$, the
volume of the active tiles is multiplied at each substitution step by a matrix conjugate to
$P(\tau^{-3})$, so the uncovered volume decays per step like the spectral radius $\rho=0.5509$ of
$P(\tau^{-3})$, a root of $x^6+732x^5-36766x^4-54904x^3+54409x^2-7588x+196$ (the norm of its
characteristic polynomial).  The two exponents are linked: $1/\beta$ is the value of $s$ at which
the spectral radius of $P(s)$ equals $1$, so $\rho<1$ is the same inequality as $\beta<\tau^3$.

\section{The golden numeration: membership without the foam}\label{sec:address}

Every point $x$ of the interior of $\B$ that is not on a cut lies in exactly one cell of the seed's
rule; mapped into the frame of that child it lies in the union of cells of the child's rule, and so
on.  The sequence of cells met is the \emph{address} of $x$, a numeration in which the digits are
cells and the places are the scales $\tau^{-k}$.  The address of a point of the interior of a
bubble ends at that bubble.

\begin{lstlisting}
procedure locate(x, K):                    # x in M, exact; K = largest scale index wanted
  tile = (D0, k=-3, id, 0)                 # the seed
  until tile is a bubble other than the seed, or its scale index exceeds K:
    (c, k, g, u) = tile; (m, t, cells) = rule(c)
    y = tau^(k+m) g^{-1}(x - u) + t        # into the frame  tau^m rep_c + t = U cells
    sigma = ( sgn(e_i . y - 1) )_{i=1..12}  # twelve exact signs in Z[tau]; a zero counts as -
    (c', g') = the cell of rule(c) with sign vector sigma, or return "none" if there is none
    tile = (c', k+m, g g', u - tau^-(k+m) A_g t)
  return tile
\end{lstlisting}

\begin{proposition}\label{prop:address}
$x\in\M$ is the centre of a bubble of scale $\tau^{-k}$ iff \texttt{locate}$(x,k)$ returns a core
of scale index $k$ and position $x$.  The descent takes at most $k+3$ steps, each consisting of
twelve sign evaluations in $\Zt$.
\end{proposition}

\begin{proof}
Every tile visited by \texttt{locate} is a tile of the foam, because the descent only passes from a
tile to one of its children \eqref{eq:children}.  So if a core $(\Dz,k,g,x)$ is returned, it is a
bubble of the foam, with centre $x$ by Lemma~\ref{lem:inM}.  Conversely let $x$ be the centre of a
bubble $T$.  Then $x$ is an interior point of $T$, hence of every ancestor of $T$, and in the frame
of the rule of an ancestor it is an interior point of exactly one cell, the one that contains $T$;
no sign is zero, the sign vector is that of this cell, and the descent follows the ancestors of $T$
and returns $T$.  The first step leads from the seed to scale index $0$, at most one step keeps the
scale index (the rule of $C_0$, with $m=0$, whose pieces are of other classes), and every other
step increases it by $m_c\ge1$; hence at most $k+2$ steps reach scale index $k$.
\end{proof}

The program checks $2\,000$ random bubble centres of the depth-$4$ foam (all recognised, $6.5$~s)
and $500$ of them asked at the scale index $k+1$ (none accepted).

\paragraph{The Galois cone.}  Where can bubble centres be?  A bubble of scale $\tau^{-k}$ with
centre $x$ satisfies $x+\tau^{-k}\Dz\subset\B$.  The arithmetic constraint is on the conjugate.
For a tile $(c,k,g,x)$ put
\[
u=(\tau^kx)^\ast=(-\tau)^{-k}x^\ast\in\R^3,
\]
the star image of the position rescaled to unit size.  Conjugating \eqref{eq:children} gives, for
every child,
\begin{equation}\label{eq:u}
u_{\text{child}}=(-\tau)^{-m_c}\,u_{\text{parent}}-R_g^\ast\,t_c^\ast .
\end{equation}
The earlier history of a tile is damped by the factors $\tau^{-m_c}$, and the last translations
are bounded, so $u$ stays bounded.

\begin{proposition}[the cone bound]\label{prop:cone}
Every tile of the foam has $|u|\le\sqrt3\,\tau^2=4.5346$, and every bubble has $|u|\le2$.  The
value $2$ is attained by the $30$ bubbles of scale $1$ inside the chocks, $x=(0,0,2)$ up to
symmetry.
\end{proposition}

\begin{proof}
Let $R_0=\sqrt3\,\tau^2$.  The tiles of scale index $0$ produced by the seed have $u=0$.  By
\eqref{eq:u} and the orthogonality of $R_g^\ast$, $|u_{\text{child}}|\le\tau^{-m_c}
|u_{\text{parent}}|+|t_c^\ast|$, and the norms $|t_c^\ast|$ are $2$ for $C_0$ (which occurs only
with $u_{\text{parent}}=0$), $\sqrt3$ for $S_0$ ($m=1$), $0.7265$ for $H_1$ ($m=1$), $1.4020$ for
$W_1$ and $1.1756$ for $H_0$ ($m=2$), $1.2361$ for $W_0$ ($m=3$).  In each case
$\tau^{-m_c}R_0+|t_c^\ast|\le R_0$, with equality for $S_0$, which proves the first bound by
induction.

For the second bound, call \emph{class tile} of $c$ the tile $(c,0,\mathrm{id},0)$; the tiles
produced by the seed are the images $(c,0,g,0)$ of the class tiles.  Iterating \eqref{eq:u}, a
bubble $b$ that descends from a tile $a$ has
\[
u_b=(-\tau)^{-(k_b-k_a)}\,u_a+R_{g_a}^\ast w,
\]
where $w$ is the value of $u$ at the corresponding bubble among the descendants of the class tile of
$c_a$, a bubble of scale index $k_b-k_a$.  Fix $K=10$.  If $k_b<K$, take for $a$ the tile of scale
index $0$ above $b$: then $|u_b|=|w|$.  If $k_b\ge K$, take for $a$ the last ancestor of $b$ with
$k_b-k_a\ge K$; since a step increases the scale index by at most $3$, $K\le k_b-k_a\le K+2$, and
$|u_b|\le\tau^{-K}R_0+|w|$.  The program \texttt{cone\_bound.py} enumerates all bubbles of scale
index at most $12$ descending from the six class tiles.  The largest $|w|$ is $2$ for the indices
below $10$ (attained only at index $0$, by the core of the chock) and $1.9015$ for the indices $10$,
$11$, $12$.  As $1.9015+\tau^{-10}R_0=1.9384<2$, the bound follows.  The norms are evaluated in
floating point, with a margin of $0.06$.
\end{proof}

So a bubble of scale $\tau^{-k}$ has $|x^\ast|\le2\tau^k$: its lift lies in the region
$\{(x,x^\ast):x\in\B,\ |x^\ast|\le2\tau^k\}$ of the six-dimensional lattice, which we call the
Galois cone because its width in internal space grows as the bubbles shrink.  For each $k$ it
contains finitely many points, and \texttt{locate} decides each one.

\begin{result}[the sieve]\label{res:sieve}
For $k=0,\dots,5$ the set $\{x\in\M: x+\tau^{-k}\Dz\subset\B,\ |x^\ast|\le2\tau^k\}$
contains \sieveca, \sievecb, \sievecc, \sievecd, \sievece, \sievecf{} points, of which
\sieveba, \sievebb, \sievebc, \sievebd, \sievebe, \sievebf{} are bubble centres by
Proposition~\ref{prop:address}.  These numbers are the exact counts $b_0,\dots,b_5$ of
Theorem~\ref{thm:beta}, as they must be by Proposition~\ref{prop:cone} (\sievet{} s in total for
the six layers).  The depth-$4$ iteration has only $1\,812$ and $4\,212$ of the $b_4=\bke$ and
$b_5=\bkf$ bubbles of layers $4$ and $5$.
\end{result}

This is the sense in which the foam can be written down by $\Zt$ arithmetic alone: the bubbles of
scale $\tau^{-k}$ are the points of $\M$ in a bounded region of the lift whose address ends at a
core of scale index $k$.  Each point is decided independently of all others, and the only recursion
is the descent along its own address.

\section{The conjugate window}\label{sec:window}

\paragraph{Definition.}  Let $\Om\subset\R^3$ be the closure of the set of values
$u=(\tau^kx)^\ast$ over all bubbles $(\Dz,k,g,x)$ of the foam, and $\Om_r$ the closure over the
bubbles with $k\equiv r\pmod 3$, $r=0,1,2$.  By Proposition~\ref{prop:cone}, $\Om$ lies in the
ball of radius $2$.  If the bubbles of scale $\tau^{-k}$ were selected among the points of $\M$ by
a window in internal space, as the vertices of a cut-and-project tiling are, that window would
have to contain all these values; this is why we call $\Om$ the \emph{window} of the bubbles.  The
recursion \eqref{eq:u} does not involve $k$, so the residues are introduced only to test whether
the window depends on the scale.  By \eqref{eq:u} the limit points of $\Om$ form the limit set of
a graph-directed system of similarities \cite{MauldinWilliams88} with ratios $\tau^{-1}$,
$\tau^{-2}$, $\tau^{-3}$, whose vertices are the pairs (class, group element) for the classes
$S_0,H_0,W_0$ and whose translations are the vectors $-R_g^\ast t_c^\ast$; the maps are composed in
the order opposite to the substitution, the most recent step being applied last.

\begin{figure}[p]\centering
\fig[0.49\linewidth]{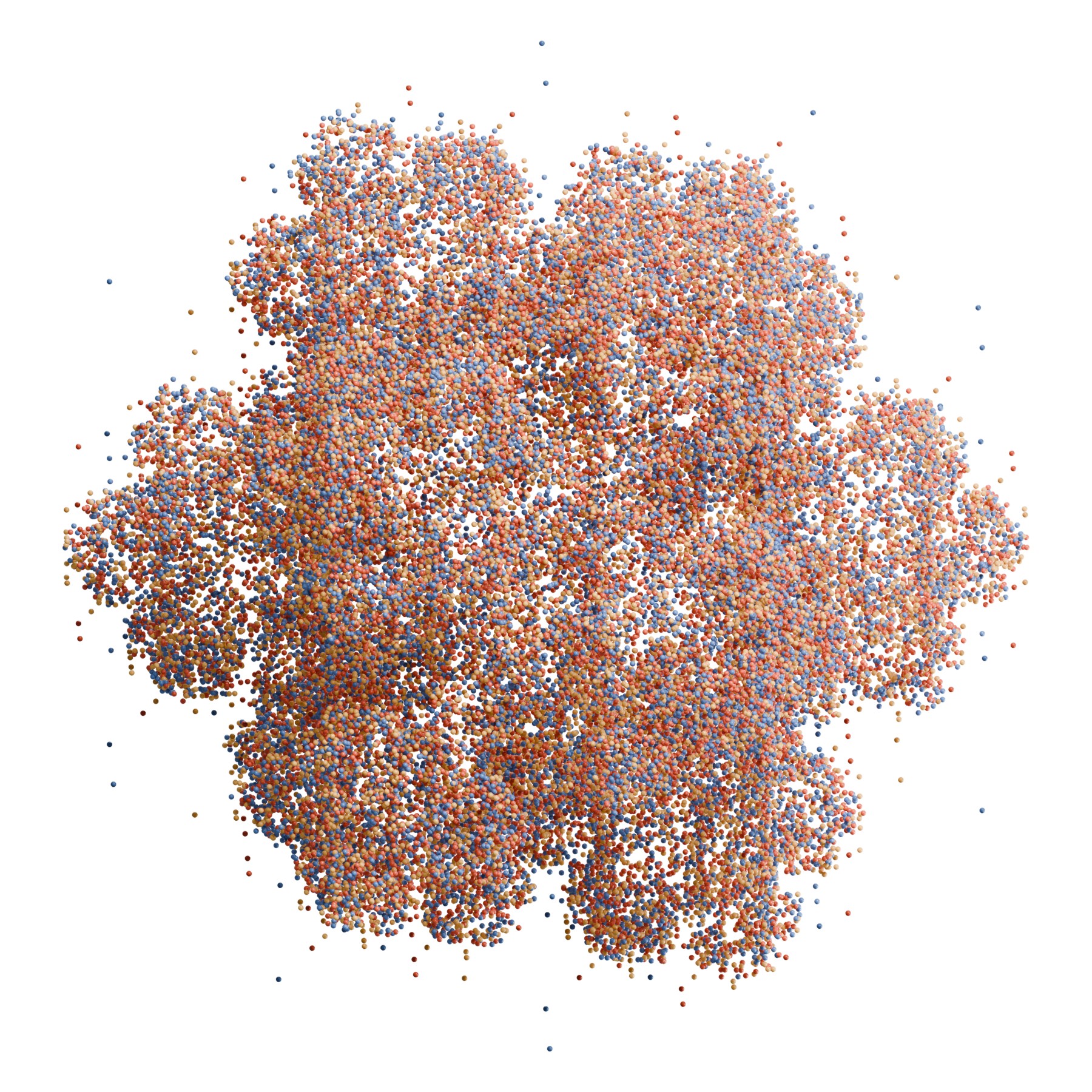}\hfill
\fig[0.49\linewidth]{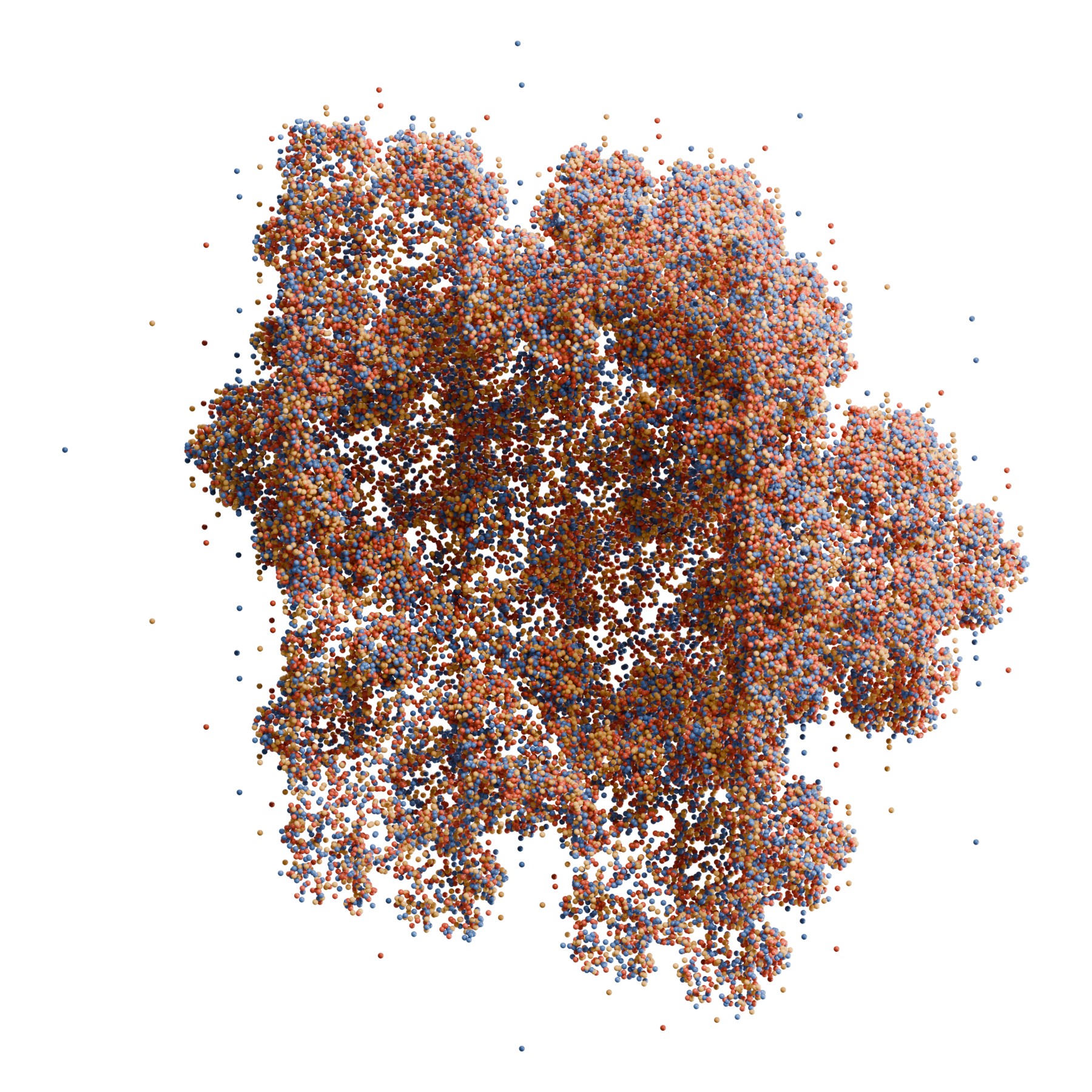}\\[4pt]
\fig[0.49\linewidth]{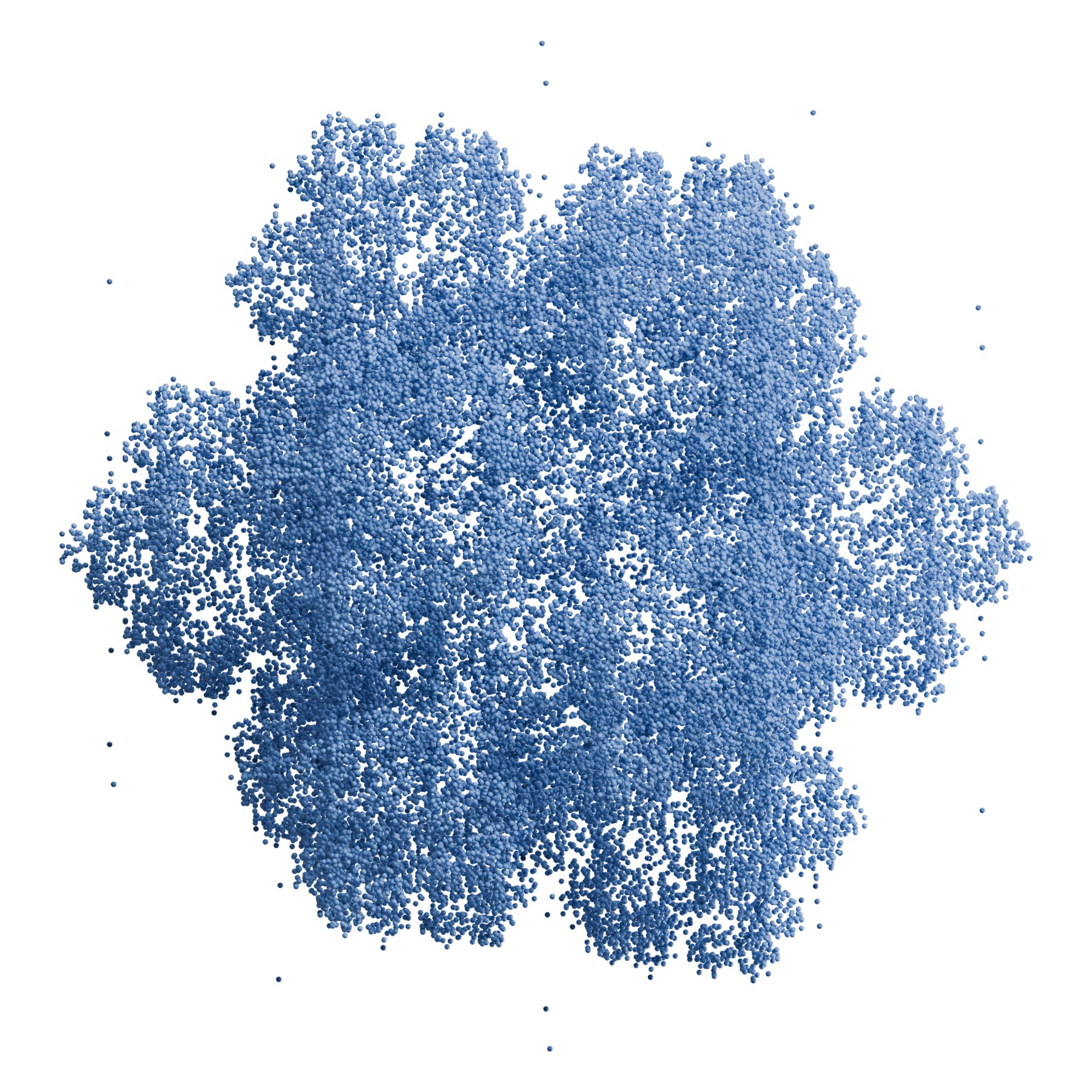}\hfill
\fig[0.49\linewidth]{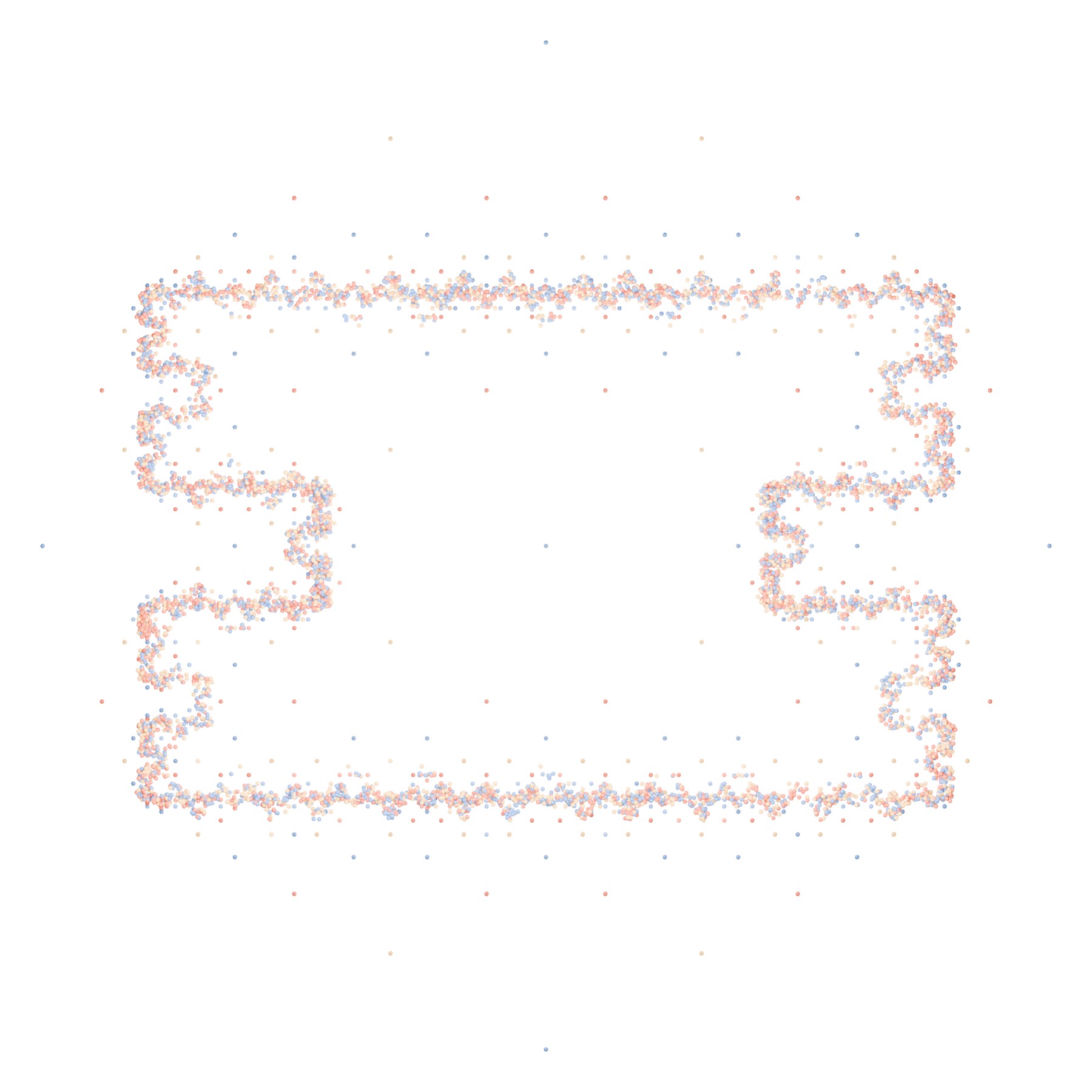}
\caption{The conjugate window $\Om$ of the Dodecafoam bubbles in internal space: $70\,000$ of the
  \nsamples{} sampled values $u=(\tau^kx)^\ast$, coloured by $k\bmod3$ (blue $0$, yellow $1$,
  red $2$).  Top left: all residues.  Top right: cut in half.  Bottom left: residue $0$ alone.
  Bottom right: the slab $|u_3|\le0.04$ around a mirror plane, seen from above; the section looks
  like a closed curve, and the isolated dots are samples lying exactly on the plane.  No sample
  has $0<|u|<0.66$.  The symmetry is that of the conjugate representation of $I_h$.}\label{fig:window}
\end{figure}

\paragraph{Sampling.}  We sample $\Om$ by \npaths{} random descents from the seed of $26$ steps
each (children chosen with probability proportional to volume), recording $u$ for every bubble
met: \nsamples{} samples, about $170\,000$ per residue, all computed exactly in $\M$ before
conjugation.  The observations (\texttt{window\_stats.py}):
\begin{itemize}
\item \emph{Radius.}  $\max|u|=2.000$ for residue $0$ (attained at the scale-$1$ chock cores),
  $1.895$ and $1.898$ for residues $1$ and $2$, in agreement with Proposition~\ref{prop:cone}.
\item \emph{One window.}  At box size $0.05$ and with samples of $80\,000$ points, the boxes hit by
  residue $0$ are hit by residue $1$ in \winovb\% of the cases and by residue $2$ in \winovc\%; a
  second, disjoint sample of residue $0$ hits \winova\%.  The three residues are the same set to
  within sampling.
\item \emph{Shape.}  The window is a shell: apart from the value $u=0$ of the central core, no
  sample has $|u|<0.66$, and the distribution of $|u|$ is $0.7\%$, $12\%$, $21\%$, $17\%$, $27\%$,
  $14\%$, $4\%$ in the bins $[0.6,0.8),\dots,[1.8,2.0]$.  Inside the shell the set is a sponge of
  lobes (Figure~\ref{fig:window}, top right).
\item \emph{Box counting.}  The samples meet \boxa, \boxb, \boxc, \boxd, \boxe{} boxes of side
  $h=0.4,0.2,0.1,0.05,0.025$.  The slopes between consecutive sizes are $2.2$ to $2.35$ until
  under-sampling sets in at the smallest size, and the covered volume $Nh^3$ decreases steadily
  ($34.6,20.1,12.8,7.6,3.3$).
\end{itemize}

\begin{theorem}\label{thm:dim}
\begin{enumerate}
\item The upper box dimension of $\Om$ satisfies $\dimB\,\Om\le\log_\tau\beta=\dimwin$, with
  $\beta=(1+\sqrt{29})/2$.  Hence $\dimH\Om\le\dimwin<3$ and $\Om$ has Lebesgue measure zero.
\item Let $X_k$ be the set of the points $\tau^kx$, $x$ the centre of a bubble of scale
  $\tau^{-k}$; these are the centres of the unit bubbles of the foam rescaled by $\tau^k$.  There
  is no relatively dense set $\Lambda\subset\R^3$ such that $X_k=\Lambda\cap\tau^k\B_k$ for all
  large $k$, where $\B_k=\{x:x+\tau^{-k}\Dz\subset\B\}$.  In particular there is no window $W$
  with nonempty interior such that the centres of the bubbles of scale $\tau^{-k}$ are, for all
  large $k$, the points $x\in\M\cap\B_k$ with $(\tau^kx)^\ast\in W$.
\end{enumerate}
\end{theorem}

\begin{proof}
(1) We use the notation of the proof of Proposition~\ref{prop:cone}, with an arbitrary $K$.  The
bubbles of scale index below $K$ are finitely many.  A bubble $b$ with $k_b\ge K$ has
$|u_b-R_{g_a}^\ast w|\le\tau^{-K}R_0$, where $w$ is the value of $u$ at one of the bubbles of scale
index $K$, $K+1$ or $K+2$ descending from one of the six class tiles, and $g_a$ is one of $120$
group elements.  The number of bubbles of scale index $j$ descending from a class tile is at most
$b_j$, because the class tiles are among the tiles of scale index $0$ of the foam; by
Theorem~\ref{thm:beta} it is at most $C\beta^j$.  So $\Om$ is covered by a finite set and at most
$C'\beta^K$ balls of radius $r=\tau^{-K}R_0$, that is, by $C''r^{-\log_\tau\beta}$ balls of radius
$r$, for every $K$.  This bounds the upper box dimension, the Hausdorff dimension is not larger,
and a set of Hausdorff dimension less than $3$ has measure zero.

(2) $\B_k$ is a dodecahedron concentric with $\B$ that tends to $\B$, so $\tau^k\B_k$ contains a
ball of radius proportional to $\tau^k$, and a relatively dense set has at least $c\,\tau^{3k}$
points in it, with $c>0$.  But $X_k$ has $b_k=O(\beta^k)$ points and $\beta<\tau^3$
(Theorem~\ref{thm:beta}).  For the last statement, $(\tau^kx)^\ast\in W$ with $x\in\M$ means
$\tau^kx\in\Lambda=\{y\in\M:y^\ast\in W\}$; we may assume $W$ bounded, since $\Om$ is, and then
$\Lambda$ is a model set, which is relatively dense \cite[Section~3]{Moody00}.
\end{proof}

The number $\log_\tau\beta$ is the similarity dimension of the graph-directed system: it is the
$s$ for which the matrix $P(\tau^{-s})$ of Section~\ref{sec:subst} has spectral radius $1$.
Equality $\dimH\Om=\log_\tau\beta$ would follow from an open set condition
\cite{MauldinWilliams88}, which we have not established; the box-counting slopes of about $2.3$
from half a million points are compatible with it, finite samples giving too small a value.

\paragraph{The model-set test.}  Although $\Om$ is not a window with interior, one may still ask
whether the condition $u(x)\in\Om$, which holds for every bubble by definition, characterises the
bubbles: does the lift of $\M$ ever meet $\Om$ at a point that is not a bubble?
Table~\ref{tab:test} compares, for the complete layers $k=2$ and $k=3$ of the depth-$4$ foam, the
distance from $u(x)$ to the nearest window sample for the bubble centres and for all other points
$x\in\M$ with $x+\tau^{-k}\Dz\subset\B$ and $|u(x)|\le2$ (outside the central core, whose interior
is not subdivided).  Non-bubbles land within $0.01$ of a sample in $5$--$12\%$ of the cases,
against $1.7\%$ for uniformly random points of the ball, and the smallest non-bubble distances
(\mindtwo{} for $k=2$, \mindthree{} for $k=3$) are below the spacing of the samples ($0.011$ to
$0.03$ between nearest neighbours).  So at the resolution of the sample the window does not
separate bubbles from non-bubbles.  We do not expect it to at any resolution, but this is not
proved (Section~\ref{sec:discussion}).  The exact criterion remains the address.

\begin{table}[ht]\centering\small
\caption{Model-set test.  Candidates: $x\in\M$ with $x+\tau^{-k}\Dz\subset\B$ outside the central
  core and $|u(x)|\le2$.  Distance $d$ of $u(x)$ to the nearest of the \nsamples{} window samples
  for the non-bubbles.  Last row: uniformly random points of the ball of radius $2$, for
  comparison.}\label{tab:test}
\makebox[\linewidth][c]{\begin{tabular}{lrrrrrr}
\toprule
 & candidates & bubbles & \multicolumn{4}{c}{non-bubbles within $d\le$} \\
 & & & $0.01$ & $0.02$ & $0.04$ & $0.08$ \\
\midrule
$k=2$ & \candtwo & \nbubtwo & \tntwoa{} (\tftwoa\%) & \tntwob{} (\tftwob\%) & \tntwoc{} (\tftwoc\%) & \tntwod{} (\tftwod\%) \\
$k=3$ & \candthree & \nbubthree & \tnthreea{} (\tfthreea\%) & \tnthreeb{} (\tfthreeb\%) & \tnthreec{} (\tfthreec\%) & \tnthreed{} (\tfthreed\%) \\
random points of the ball & & & $1.7\%$ & $9.3\%$ & $24.6\%$ & $41.0\%$ \\
\bottomrule
\end{tabular}}
\end{table}

\section{Comparison with the construction of 1999}\label{sec:compare}

Goodman-Strauss's article \cite{GS99} builds the Dodecafoam in three stages.  The cells come from
the twelve planes faceting the dodecahedron $D_3$, our $\B$ (his Fig.~4).  The rules are stated in
an appendix as lists of isometries, written as words in a few generating rotations and
reflections.  The substitution is iterated with a Mathematica routine like his \texttt{SubIt},
which keeps a list of objects $\{\text{position},\text{type}\}$, the position being a $4\times4$
real matrix in homogeneous coordinates, and applies to each object its rule, one matrix product
per child; the published images were computed by a faster program of Dan Krech at the Geometry
Center and viewed in Geomview.  How the rules were found is stated plainly: ``I had been making
hand-drawn (!) pictures and paper models of dodecafoam for a few years'', and ``it is unclear to me,
frankly, whether the particular choices made are art or science''.  Our geometric search replaced
the hand by an exhaustive search in exact coordinates and reproduced his rules.
Table~\ref{tab:compare} puts the three constructions side by side, and the rest of this section
measures the entries that are quantitative.

\begin{table}[t]\centering\small
\begin{tabular}{@{}p{2.2cm}p{4.0cm}p{3.5cm}p{4.9cm}@{}}\toprule
 & Goodman-Strauss 1999 \cite{GS99} & geometric search & this article\\\midrule
cells & $165$ cells of twelve planes & sign vectors, exact geometry & sign vectors over twelve
  cuts; the group permutes the signs\\
rules & by hand, over several years; one rule per class chosen & exhaustive geometric search, minutes;
  $22$ rules & solutions of $\lambda c+e\cdot t\in\{\pm1,\pm\tau^3\}$ in $\B$
  (Theorem~\ref{thm:rules}), seconds; complete for all $\lambda>0$ (Result~\ref{res:rules})\\
isometries & words in generators, as $4\times4$ real matrices & $3\times3$ matrices
  over $\tfrac12\Zt$ & permutations of the twelve cuts, matrices in $\GL(3,\Zt)$\\
a tile & $12$ reals and a type & $12$ reals and a class & six integers, a group index, a class\\
iteration & one real $4\times4$ product per child & C++ in double precision & the
  same recursion in integers; same cost, exact (Fig.~\ref{fig:perf}a,d)\\
membership & --- & --- & the address: at most $k+3$ exact sign vectors for scale $\tau^{-k}$
  (Fig.~\ref{fig:perf}c)\\
a layer & by iterating everything & same & the sieve over the Galois cone, without the foam
  (Fig.~\ref{fig:perf}b)\\
counting & --- & bubble volume per step & $b_k$ exactly; growth rate $\beta=(1+\sqrt{29})/2$\\
window & --- & --- & $\dimB\,\Om\le\log_\tau\beta$, measure zero, no window with interior\\
other foams & seeds of Figs.~7--11 & --- & any union of cells as seed, any rule choice from the
  list (Section~\ref{sec:other})\\
\bottomrule
\end{tabular}
\caption{The three constructions of the Dodecafoam.}\label{tab:compare}
\end{table}

\paragraph{What is faster and what is not.}  We implemented an equivalent of \texttt{SubIt}: a tile
is a $4\times4$ real matrix $T$, a rule is a fixed list of matrices
$K_i=\bigl(\begin{smallmatrix}R_i/\lambda&0\\-t/\lambda&1\end{smallmatrix}\bigr)$, a child is
$K_iT$.  It uses the same rules and the same group elements $R_i$ as the integer version, and both
run in Python on one core (\texttt{tools/compare.py}).  Figure~\ref{fig:perf} shows the result.  The
substitution itself costs about the same in integers: a child costs one $4\times4$ real product in one
method and one $6\times6$ integer product plus a permutation lookup in the other, about $\puge$
against $\puze\,\mu$s per tile at depth~$4$ (the integer version is some $15\%$ faster), and both scale linearly with the number of tiles
($\pitf$ tiles at depth~$5$ in $\ptgf$ and $\ptzf$ seconds).  The floating-point drift of the real
method is $\perrf\times10^{-15}$ at depth~$5$, harmless for a picture.  What the integers give is
exact equality: two tiles coincide iff their six integers, group index and class coincide, and no
tolerance enters.  The identities of this article (the $22$ rules, the equality of the bubble sets
with the geometric iteration, the sub-foam property of Section~\ref{sec:other}) are equalities of
integers.

The gain is in the two operations that the 1999 construction does not have.  A membership query
costs at most $k+3$ exact sign vectors: $\pada$ to $\padh$ ms for bubble centres of scale index
$k=0,\dots,7$, growing linearly in $k$ and independent of the size of the foam
(Fig.~\ref{fig:perf}c); the iteration answers the same question only after producing everything
down to that scale.  For complete layers the comparison is between two growth rates.  The sieve
examines the points of the Galois cone, whose number grows by $\tau^3$ per layer, each decided in
at most $k+3$ steps.  The iteration needs depth $k+1$ to complete layer $k\ge1$ ($31,32,212,392,
1832$ bubbles of scales $\tau^0,\dots,\tau^{-4}$ are complete at depths $2,2,3,4,5$), and the
number of tiles grows by about $16$ per depth.  In our implementation the sieve is the slower
method for the layers we ran: \sievet\ s for layers $0$ to $5$, against $\ptzf$ s for the iteration
to depth~$5$, which contains layers $0$ to $4$ and all but $20$ of the $4\,532$ bubbles of
layer~$5$.  The cost of the last sieved layer was $\times\sieveratio$ that of the one before.
Extrapolating the sieve by its worst-case growth $\tau^3(k+1)/k$ per layer and the iteration by
$16$ per depth, the sieve becomes the faster method from layer $\crossk$ on (Fig.~\ref{fig:perf}b);
layers beyond $5$ were not run.  Both the sieve and the address are local: they can be run on a
region of space, or on a single point, whereas the iteration must produce at least all the
ancestors of the tiles of that region.

\paragraph{What is simpler.}  The input of the whole construction is the Gram matrix of the twelve
normals and the two levels $1$ and $\tau^3$; the Euclidean geometry survives only as the inner
product $e\cdot x$ on $\M$, evaluated exactly in $\Zt$.  Every other object is combinatorial or
integral: a cell is a sign vector, an isometry a permutation of twelve signs, a rule a solution of
a linear system with a finite right-hand side, a tile a six-tuple, a layer a set of six-tuples cut
out by a bounded region of the lift and a descent.  The choice of rules, which was ``art or
science'' in 1999, becomes a choice from a finite list that is complete
(Result~\ref{res:rules}), and Section~\ref{sec:other} examines all the choices.

\pgfplotsset{compat=1.18,
  perf/.style={width=0.47\linewidth, height=0.36\linewidth, tick label style={font=\scriptsize},
    label style={font=\scriptsize}, title style={font=\scriptsize\bfseries, yshift=-2pt},
    legend style={font=\scriptsize, draw=none, fill=white, fill opacity=0.85, text opacity=1, at={(0.97,0.03)}, anchor=south east, cells={anchor=west}},
    legend cell align=left, axis line style={gray!60}, tick style={gray!60}, grid=major, grid style={gray!18},
    every axis plot/.append style={line width=1.1pt, mark size=1.6pt}}}
\definecolor{cZ}{HTML}{4F77B5}\definecolor{cGS}{HTML}{D9573B}\definecolor{cSieve}{HTML}{2F8F62}\definecolor{cAddr}{HTML}{8A6FC4}
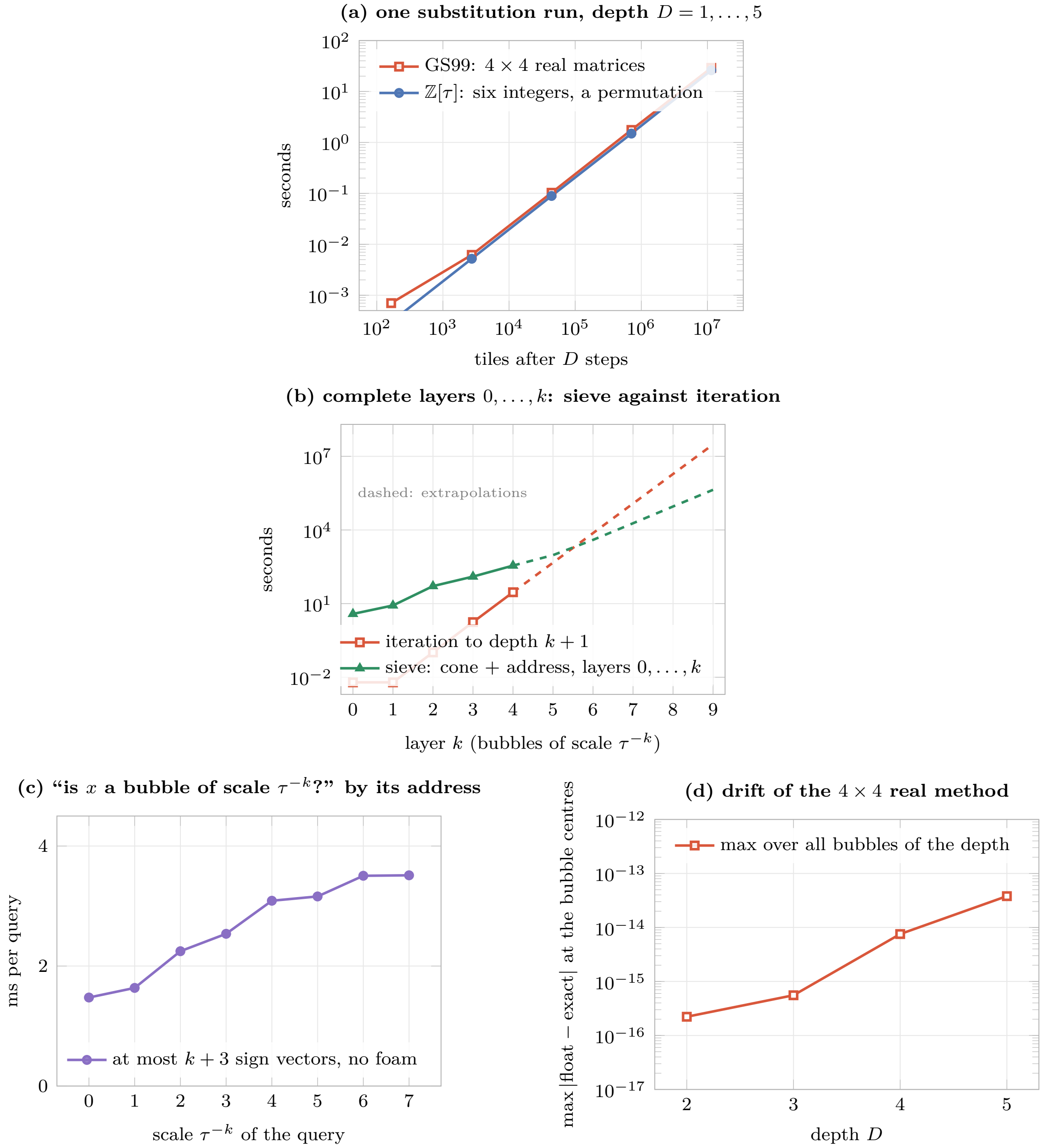
\begin{figure}[tp]\centering
\begin{tikzpicture}
\begin{axis}[perf, xmode=log, ymode=log, xlabel={tiles after $D$ steps}, ylabel={seconds},
  title={(a) one substitution run, depth $D=1,\dots,5$}, ymin=0.0005, ymax=100, legend style={at={(0.03,0.97)}, anchor=north west},
  log basis x=10, xtick={100,1000,10000,100000,1000000,10000000}, xticklabels={$10^2$,$10^3$,$10^4$,$10^5$,$10^6$,$10^7$}]
\addplot[cGS, mark=square*, mark options={fill=white}] table[x=tiles, y=tgs] {data/perf_iter.dat};
\addlegendentry{GS99: $4\times4$ real matrices}
\addplot[cZ, mark=*] table[x=tiles, y=tz] {data/perf_iter.dat};
\addlegendentry{$\Zt$: six integers, a permutation}
\end{axis}
\end{tikzpicture}\hfill
\begin{tikzpicture}
\begin{axis}[perf, ymode=log, xlabel={layer $k$ (bubbles of scale $\tau^{-k}$)}, ylabel={seconds},
  title={(b) complete layers $0,\dots,k$: sieve against iteration}, xmin=-0.3, xmax=9.3, ymin=0.002, ymax=2e8, xtick={0,...,9}]
\addplot[cGS, dashed, mark=none, forget plot] table[x=k, y=titer] {data/perf_extrap.dat};
\addplot[cSieve, dashed, mark=none, forget plot] table[x=k, y=tsieve] {data/perf_extrap.dat};
\addplot[cGS, mark=square*, mark options={fill=white}] table[x=k, y=tgs] {data/perf_layers.dat};
\addlegendentry{iteration to depth $k+1$}
\addplot[cSieve, mark=triangle*] table[x=k, y=tsieve] {data/perf_layers.dat};
\addlegendentry{sieve: cone $+$ address, layers $0,\dots,k$}
\node[font=\tiny, gray, anchor=west] at (axis cs:-0.1,3e5) {dashed: extrapolations};
\end{axis}
\end{tikzpicture}\\[2pt]
\begin{tikzpicture}
\begin{axis}[perf, xlabel={scale $\tau^{-k}$ of the query}, ylabel={ms per query}, ymin=0, ymax=4.5, xtick={0,...,7},
  title={(c) ``is $x$ a bubble of scale $\tau^{-k}$?'' by its address}]
\addplot[cAddr, mark=*] table[x=k, y=ms] {data/perf_address.dat};
\addlegendentry{at most $k+3$ sign vectors, no foam}
\end{axis}
\end{tikzpicture}\hfill
\begin{tikzpicture}
\begin{axis}[perf, ymode=log, xlabel={depth $D$}, ylabel={$\max|\text{float}-\text{exact}|$ at the bubble centres}, xtick={1,...,5},
  ymin=1e-17, ymax=1e-12, title={(d) drift of the $4\times4$ real method}, legend style={at={(0.03,0.97)}, anchor=north west}]
\addplot[cGS, mark=square*, mark options={fill=white}] table[x=depth, y=maxerr] {data/perf_error.dat};
\addlegendentry{max over all bubbles of the depth}
\end{axis}
\end{tikzpicture}
\caption{Measured on one core in pure Python, same rules, same group elements.  (a) The substitution
  itself costs the same by both methods: about $\puge$ and $\puze\,\mu$s per tile at depth~$4$.  (b) All
  layers up to $k$: the iteration to depth $k+1$ ($\times16$ tiles per depth) against the sieve over the
  Galois cone ($\times\tau^3$ candidates per layer, each decided in at most $k+3$ steps).  The sieve is
  slower for the layers measured ($k\le5$); the dashed lines extrapolate the iteration by $\times16$ and
  the sieve by its worst-case growth $\times\tau^3(k+1)/k$ per layer, and cross before layer~$\crossk$.
  (c) The cost of one membership query for a bubble centre grows linearly in $k$ (one exact sign vector
  per step).  (d) The $\Zt$ version is exact; the real version drifts by at most $4\cdot10^{-14}$ at
  depth~$5$, which is harmless for pictures and useless for identities.}\label{fig:perf}
\end{figure}

\section{Other foams}\label{sec:other}

The 1999 article shows, besides the Dodecafoam filling $D_3$, the same substitution filling a hat,
a wedge, a spike and a chock (its Fig.~7), the small stellated dodecahedron (Fig.~8), the great
dodecahedron (Fig.~10) and the great stellated dodecahedron (Fig.~11), and remarks that ``these
solids are nested, exactly as in Fig.~5''.  In the present setting two parameters of the
construction are choices, and both are data for the same arithmetic: the seed, and the rule
selected for each class from the complete list of Result~\ref{res:rules}.

\begin{figure}[p]\centering
\fig[0.49\linewidth]{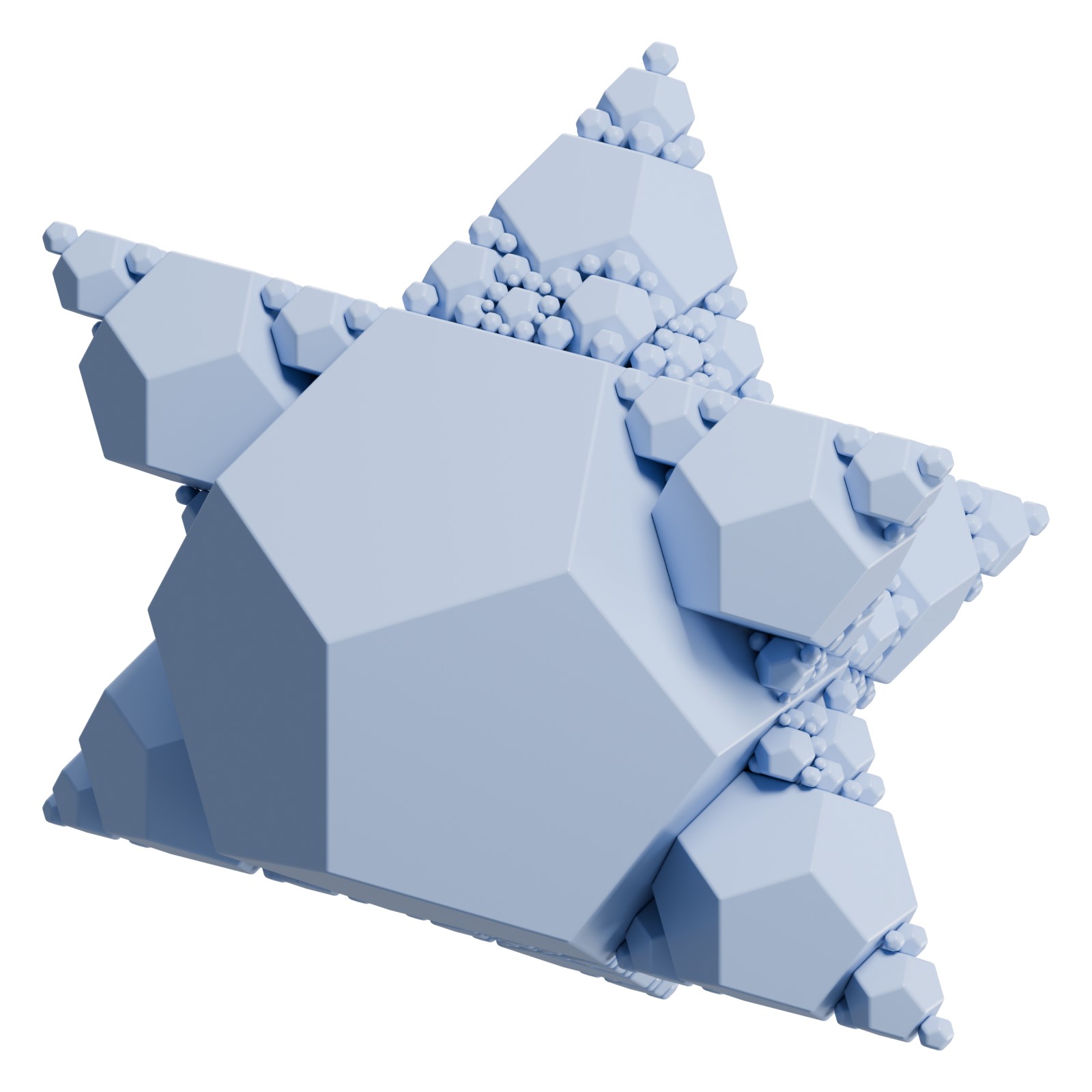}\hfill
\fig[0.49\linewidth]{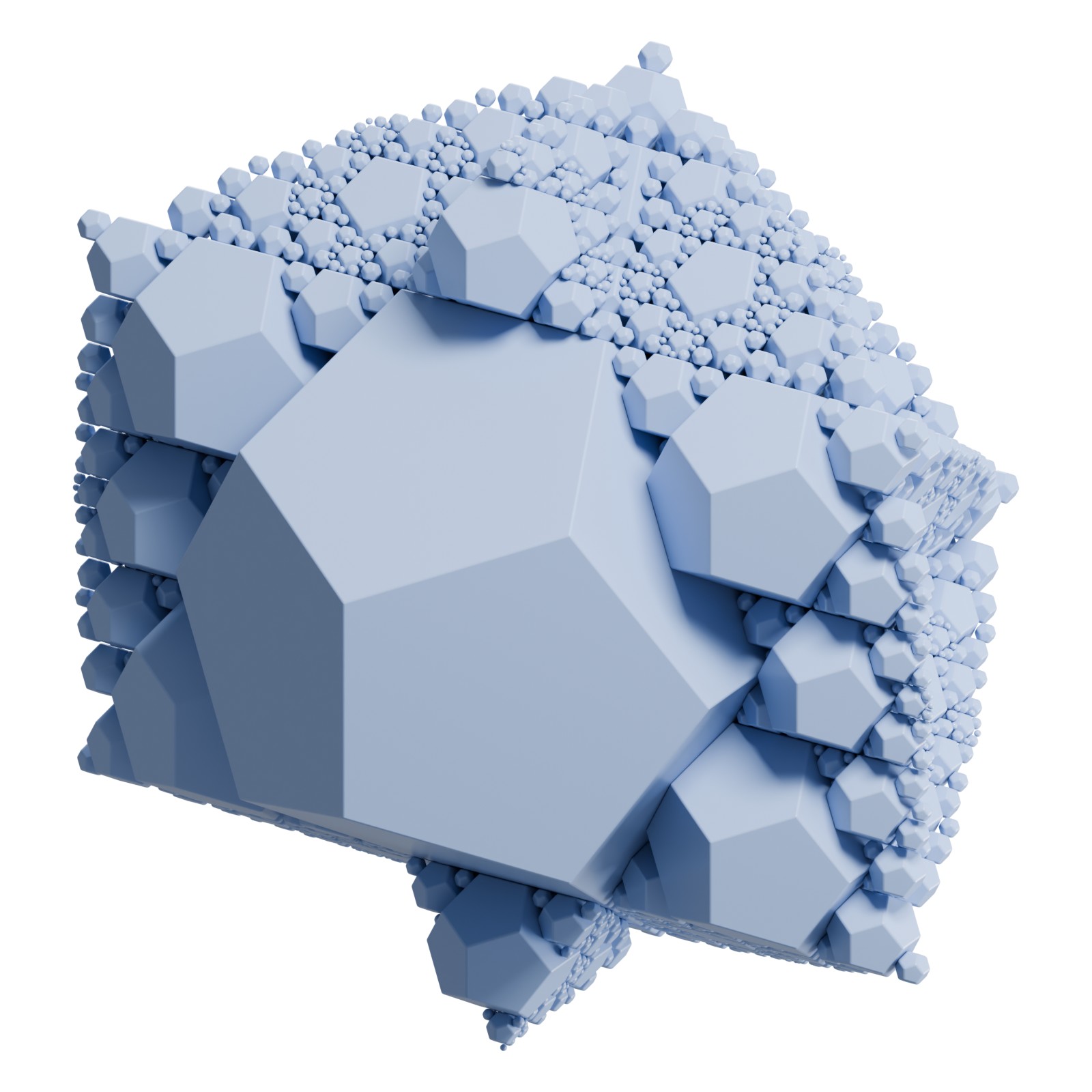}\\[4pt]
\fig[0.49\linewidth]{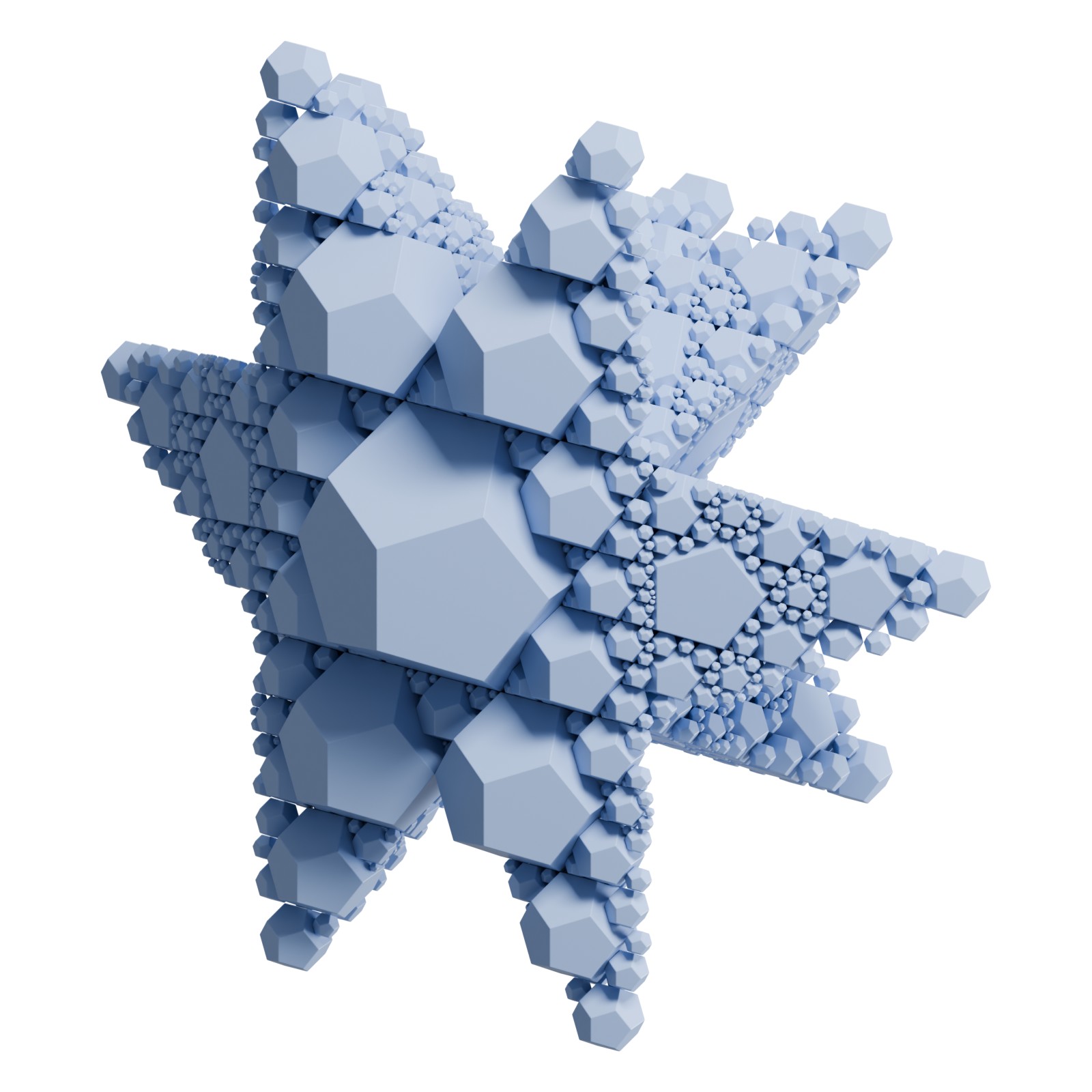}\hfill
\fig[0.49\linewidth]{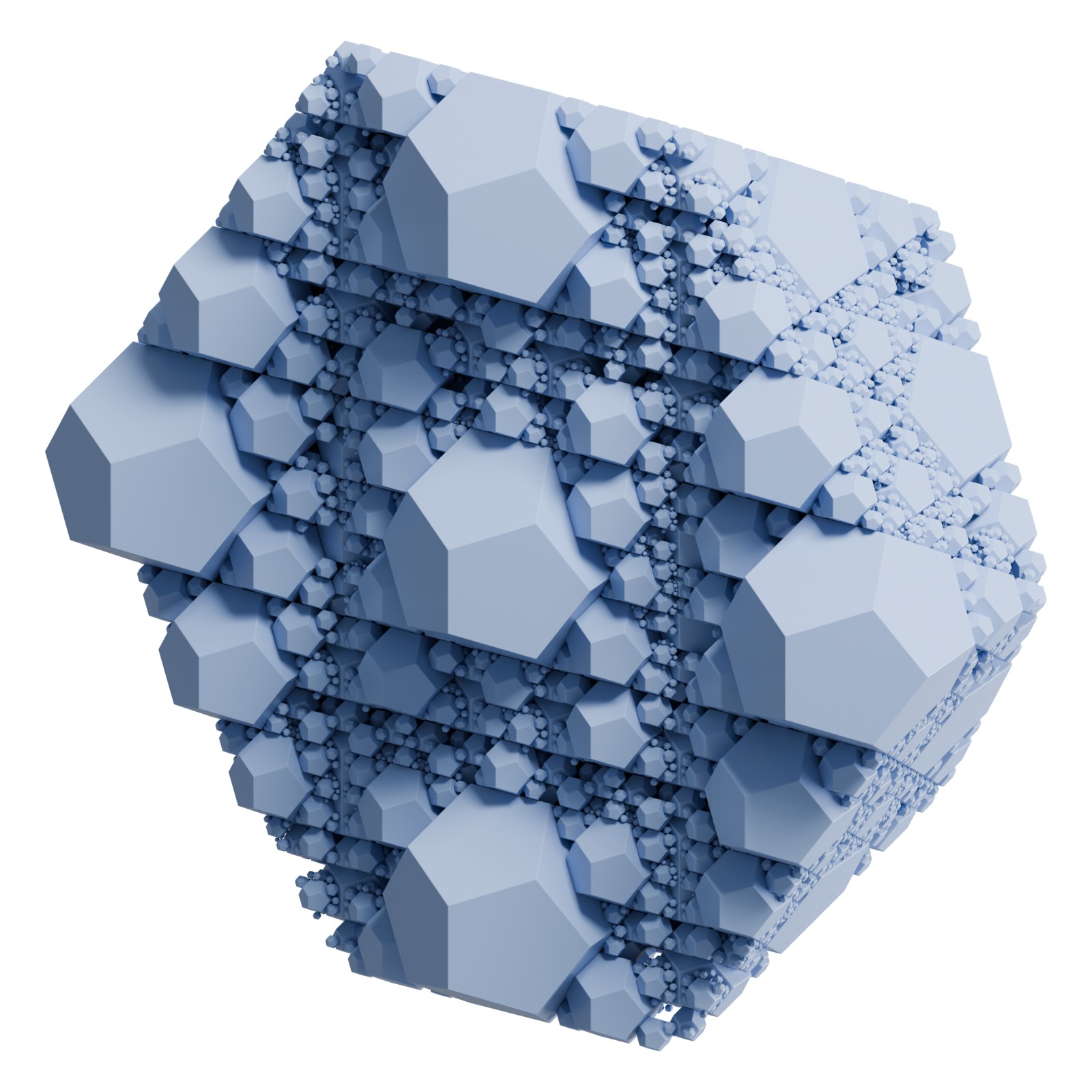}
\caption{Other foams by the same integer arithmetic, depth $4$, bubbles only, cut in half.  Top:
  the small stellated dodecahedron ($\Dz+12H_0$, \ofssdbub{} bubbles) and the great dodecahedron
  ($+30W_0$, \ofgdbub{}) filled with dodecafoam---Figures~8 and 10 of \cite{GS99}---as sub-foams of
  the Dodecafoam.  Bottom left: the great stellated dodecahedron ($+20S_0$, \ofgsdbub{}), his
  Figure~11.  Bottom right: a foam that is not in the 1999 article: the spike substitution of
  Table~\ref{tab:variants} with seed $\B$, \vrspikebub{} bubbles, growth rate $\vrspikebeta$.}\label{fig:other}
\end{figure}

\paragraph{Seeds: sub-foams of one foam.}  The cells of the arrangement are the $165$ sign
vectors; the first substitution step turns the seed $\tau^3\Dz$ into exactly these cells at scale
$\tau^0$.  Every solid of the list above is a union of cells---$\Dz$, $\Dz+12H_0$ (small stellated),
$\Dz+12H_0+30W_0$ (great), $\Dz+12H_0+30W_0+20S_0$ (great stellated), all $165$ (the dodecahedron
$D_3$)---so its foam is the set of tiles descending from a subset of the sign vectors, and it is a
sub-foam of the Dodecafoam in the exact sense: the \ofgsdtiles{} tiles of the depth-$4$ foam of the
great stellated dodecahedron are tiles of the depth-$4$ Dodecafoam ($\ofdtiles$ tiles), as sets
of integers.  Table~\ref{tab:seeds} gives the counts and Figure~\ref{fig:other} the renders; all are
produced by \texttt{tools/otherfoams.py} from the integer tiles.  A single cell is a seed too: the
foams of $H_1$ and $W_1$ are those of $H_0$ and $W_0$ shifted by one scale, as the rules
$\tau H_0=H_1$, $\tau W_0=W_1$ say.

\begin{table}[t]\centering\small
\begin{tabular}{@{}lrrrrrr@{}}\toprule
seed & cells & tiles at depth $4$ & bubbles & \multicolumn{3}{c}{bubble fraction, $2,3,4$ steps}\\\midrule
dodecahedron $D_3$ (all cells) & 165 & \ofdtiles & \ofdbub & 0.565 & 0.769 & \ofdfrac\\
great stellated dodecahedron & 63 & \ofgsdtiles & \ofgsdbub & 0.443 & 0.688 & \ofgsdfrac\\
great dodecahedron & 43 & \ofgdtiles & \ofgdbub & 0.570 & 0.760 & \ofgdfrac\\
small stellated dodecahedron & 13 & \ofssdtiles & \ofssdbub & 0.746 & 0.866 & \ofssdfrac\\
one small hat $H_0$ & 1 & \ofhtiles & \ofhbub & 0.541 & 0.758 & \ofhfrac\\
one small wedge $W_0$ & 1 & \ofwtiles & \ofwbub & 0.286 & 0.589 & \ofwfrac\\
one spike $S_0$ & 1 & \ofstiles & \ofsbub & 0.403 & 0.665 & \ofsfrac\\
one chock $C_0$ & 1 & \ofctiles & \ofcbub & 0.686 & 0.848 & \ofcfrac\\
\bottomrule
\end{tabular}
\caption{The foams of Goodman-Strauss's Figures~7--11 as sub-foams of the Dodecafoam, with the
  fraction of the seed's volume filled by bubbles (rules of 1999;
  depth $4$ counted from the seed $\B$, so that the cells of a seed are at depth $1$).}\label{tab:seeds}
\end{table}

\paragraph{Rules: the same foam, or a new one.}  Result~\ref{res:rules} gives, counting
placements and leaving out the single-cell solutions at scale $1$, $3,5,6,8,7,10$ rules for the
classes $C_0,S_0,H_1,W_1,H_0,W_0$, hence $50\,400$ substitutions with one rule per class (the core
is expanded only as the seed).  The 1999 choice takes for each class the rule containing a core at
the smallest scale.  Three other choices (\texttt{tools/otherfoams.py},
Table~\ref{tab:variants}):

\begin{itemize}
\item \emph{chock}: every rule contains a chock and no core ($S_0\to C_0+H_1+3W_1$ at $\tau$,
  $H_1\to C_0+2W_1$ at $\tau$, $W_1\to C_0+2H_1+6W_1$ at $\tau^2$, $H_0\to C_0+2W_1$ at $\tau^2$,
  $W_0\to C_0+2H_1+6W_1$ at $\tau^3$); bubbles arise only through $C_0\to\Dz+4H_0+W_0$.
\item \emph{wedge}: $W_0\to W_1$, $W_1\to H_1+2W_1$, $H_1\to C_0+2W_1$ (all at $\tau$), the rest
  as in 1999.
\item \emph{spike}: $S_0\to 3C_0+\Dz+5S_0+3W_1+9H_0+15W_0$ at $\tau^2$ instead of the $14$-cell rule
  at $\tau$, the rest as in 1999.
\end{itemize}

\begin{result}\label{res:variants}
The chock and wedge substitutions produce the same bubbles as the Dodecafoam as far as tested:
their bubble counts $b_k$ coincide with those of Theorem~\ref{thm:beta} for $k\le12$, and, decided
by the address in both directions,
every one of the $55\,407$ bubbles of the chock foam at depth $6$ (every one of the
\vrwedgebub{} of the wedge foam at depth $4$) is a bubble of the Dodecafoam of the same scale, and
every one of the $43\,933$ bubbles of the depth-$4$ Dodecafoam is a bubble of each variant.  The
spike substitution produces a different foam: $b_k=31,12,272,132,2692,3072,\dots$, growth rate
$\vrspikebeta$ (against $3.1926$), $\log_\tau$ of it $\vrspikedim$; at depth~$4$
it shares $53\,353$ of its \vrspikebub{} bubbles with the Dodecafoam.
\end{result}

In the chock and wedge foams the bubbles appear later: every rule of the chock substitution
produces one chock, which opens into a core one step later, so its bubbles at depth $6$ are as
many as its tiles at depth $4$.

Some choices therefore do not change the bubbles and others do.  A rule containing a chock
re-parses the same hierarchy (the chock opens into a core at the same scale,
$C_0=\Dz+4H_0+W_0$): the depth at which the bubbles appear changes, the bubbles do not.  Replacing
the rule of the spike by the one at scale $\tau^2$ changes the bubble set and, in this variant, the
growth rate.  The quantities of Sections~\ref{sec:subst}--\ref{sec:window} are available for any
of the $50\,400$ substitutions by the same code: the delay recursion with the variant's delays and
the address with the variant's rules.

\paragraph{Listing the foams.}  Since the list of Result~\ref{res:rules} is finite, the
substitutions it allows can be listed and compared.  A substitution of the Dodecafoam kind is a
point of the product of the six solution sets: one rule (scale and placement) for each of
$C_0,S_0,H_1,W_1,H_0,W_0$, the single-cell solutions at scale $1$ excluded, the seed rule
$\tau^3\Dz=\B$ fixed and the core frozen.  There are $3\cdot5\cdot6\cdot8\cdot7\cdot10=50\,400$
such points, and $2\cdot3\cdot2\cdot5\cdot3\cdot6=1\,080$ rule shapes when placements are
identified.  Two exact invariants are attached to each point by \texttt{tools/foamlist.py}: the
count sequence $(b_k)$ from the delay recursion, and the \emph{fingerprint} $F(A,K)$, the set of
bubbles (scale index $\le K$, centre in $\M$) of the infinite foam of the representative tile of
each class $c$, by the memoised recursion
\begin{multline*}
F(c,K)=\{(m,v):\text{one for each core of the rule}\}\\
\cup\ \{(m+k',\ \tau^{-m}A_{g_i}u'+v):(k',u')\in F(c_i,K-m)\},
\end{multline*}
where $\tau^m\mathrm{rep}_c+t=\bigcup g_i(\mathrm{rep}_{c_i})$ is the rule chosen for the class
$c$ and $v=-\tau^{-m}t$.  The foam of any seed is the union
of the transformed fingerprints of its cells, so equal fingerprints for all classes mean equal
bubbles down to $K$ scale indices below every cell, and different fingerprints mean different
foams.  With $K=6$:

\begin{result}\label{res:enlist}
The $50\,400$ substitutions produce $18\,899$ distinct bubble sets ($378$ for the $1\,080$ shapes)
and $238$ distinct count sequences.  The substitutions with the fingerprint of the Dodecafoam are
exactly those with $C_0\to\Dz+4H_0+W_0$ and the spike at scale $\tau$ (either rule, any
placement), with every other choice free: a product set of
$1\cdot4\cdot6\cdot8\cdot7\cdot10=13\,440$ substitutions ($360$ shapes).  Two members of this set
that differ in every free choice, the rules of 1999 and the substitution with every other class at
its smallest scale and chock-type rules, were also compared by the address: all $43\,933$ bubbles of
the first at depth~$4$ are bubbles of the second, and all $95$ bubbles of the second at depth~$4$
are bubbles of the first.  $11\,340$ substitutions
are dead (no rule reaches a core).  Among the rest, $\beta$ ranges from $\sqrt6=2.4495$ ($\log_\tau\beta=1.862$: bubbles at even layers only, $b_{2j}=12\cdot6^{\,j-1}$ for $j\ge1$) to $4.2339$ ($\log_\tau\beta=2.9989$,
almost $\tau^3$).  $16\,849$ bubble sets come from a single substitution each: placements rotate the
sub-foam of each tile in its own frame and give different sets.
\end{result}

\begin{table}[t]\centering\footnotesize
\begin{tabular}{@{}lp{4.6cm}llp{3.4cm}r@{}}\toprule
foam & distinguishing rules & $\beta$ & $\log_\tau\beta$ & $b_0,b_1,\dots$ & shapes / subst.\\\midrule
dead & $C_0\to\tau^1$ and no core rule & --- & --- & $1,0,0,\dots$ & 40 / 11\,340\\
sparsest & C01.S00.H10.W12.H02.W01 & $\sqrt6$ & 1.862 & $1,0,12,0,72,0,432$ & 6 / 450\\
Dodecafoam & C00 and $S_0$ at $\tau$, the rest free & $\tfrac{1+\sqrt{29}}2$ & 2.412 & $31,32,212,392,1832$ & 360 / 13\,440\\
spike at $\tau^2$ & C00.S02.H11.W14.H00.W01 & 3.3001 & 2.481 & $31,12,272,132,2692$ & 4 / 6\\
densest & C01.S00.H10.W12.H02.W03 & 4.2339 & 2.9989 & $1,0,12,0,102,30,884$ & 2 / 180\\
\bottomrule
\end{tabular}
\caption{Five of the $378$ foams of the $1\,080$ rule shapes (Result~\ref{res:enlist}).  A code such
  as S02 names the rule of the class $S_0$ with index $2$ in the list of \texttt{out/foamlist\_shapes.txt},
  in which the rules of a class are ordered by scale: C00 is $C_0=\Dz+4H_0+W_0$ and C01 the rule at
  scale $\tau$; S00, S01 are the two rules of the spike at scale $\tau$ and S02 the one at $\tau^2$.
  Last column: numbers of shapes and of substitutions.  Figure~\ref{fig:betas} shows the growth
  rates of all the foams that have one.}\label{tab:enlist}
\end{table}

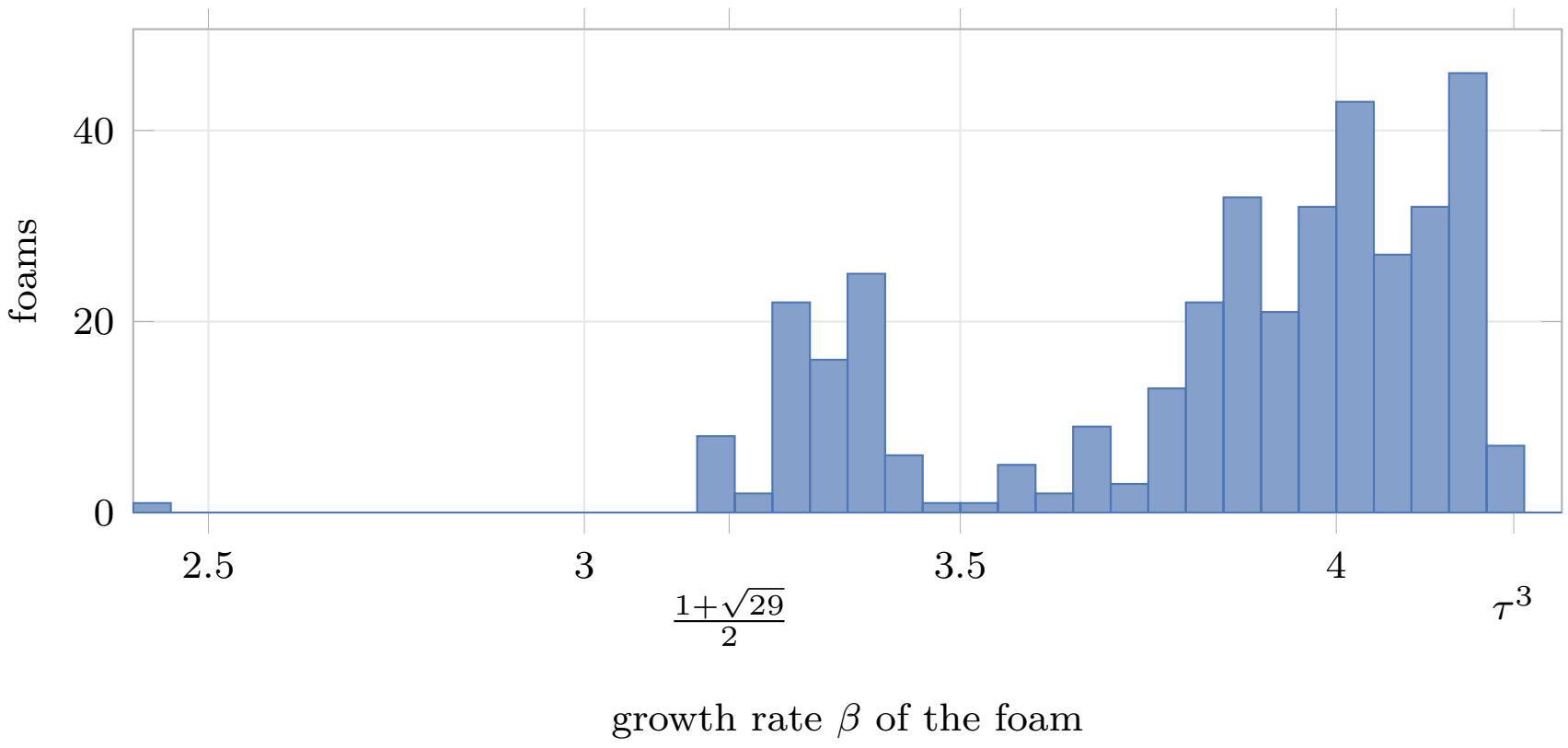
\begin{figure}[t]\centering
\begin{tikzpicture}
\begin{axis}[width=0.75\linewidth, height=0.32\linewidth, ybar interval, x tick label as interval=false, ymin=0, xmin=2.4, xmax=4.3, enlarge x limits=false,
  xlabel={growth rate $\beta$ of the foam}, ylabel={foams}, tick label style={font=\scriptsize},
  label style={font=\scriptsize}, axis line style={gray!60}, tick style={gray!60}, grid=major, grid style={gray!18},
  xtick={2.5,3,3.5,4}, extra x ticks={3.1926,4.2361}, extra x tick labels={$\tfrac{1+\sqrt{29}}2$,$\tau^3$},
  extra x tick style={tick label style={font=\scriptsize, yshift=-7pt}, grid=none}]
\addplot[fill=cZ!70, draw=cZ, hist={bins=38, data min=2.4, data max=4.3}] table[y=beta] {data/foam_betas.dat};
\end{axis}
\end{tikzpicture}
\caption{The growth rates of the foams of the $1\,080$ rule shapes (the $40$ dead shapes
  excluded).  The Dodecafoam's $\beta=(1+\sqrt{29})/2$ is marked.  All values are below $\tau^3$, as
  they must be for disjoint bubbles (Corollary~\ref{cor:density}).}\label{fig:betas}
\end{figure}

In this sense a foam is a class of points of the product of the solution sets under the
equivalence ``same fingerprint'', and for the Dodecafoam that class is a product: two choices decide
it, the other four are free (Table~\ref{tab:enlist}).  Every candidate substitution is in the
list, with its growth rate and first layers, and whether two of them differ is decided by the
fingerprint: exactly when the fingerprints differ, and down to $K$ scale indices when they agree,
with the address available for a check of any pair to a given depth.  The growth rate does not
determine the foam: there are substitutions with the spike at scale $\tau^2$ whose growth rate is
again $(1+\sqrt{29})/2$ and whose bubbles differ from those of the Dodecafoam.

\begin{table}[t]\centering\small
\begin{tabular}{@{}lllrrr@{}}\toprule
rules & $\beta$ & $\log_\tau\beta$ & tiles at depth $4$ & bubbles & same bubbles as 1999\\\midrule
1999 (core, smallest scale) & \vrgsbeta & \vrgsdim & \vrgstiles & \vrgsbub & ---\\
chock & \vrchockbeta & \vrchockdim & \vrchocktiles & \vrchockbub & yes (all bubbles tested)\\
wedge & \vrwedgebeta & \vrwedgedim & \vrwedgetiles & \vrwedgebub & yes (all bubbles tested)\\
spike at $\tau^2$ & \vrspikebeta & \vrspikedim & \vrspiketiles & \vrspikebub & no\\
\bottomrule
\end{tabular}
\caption{Four substitutions from the list of Result~\ref{res:rules}, seed $\tau^3\Dz$
  (Result~\ref{res:variants}).}\label{tab:variants}
\end{table}

\section{Four dimensions}\label{sec:4d}

The equivalence of Theorem~\ref{thm:rules} and its proof do not depend on the dimension: they
use only that the cells are the pieces into which finitely many cutting hyperplanes divide a convex
polytope.  The four-dimensional analogue of the stellation dissection is the \emph{star cut}
\cite{SDR}: a $120$-cell $\B$ with facet hyperplanes $n\cdot x=2\tau^2$ is cut by the $120$
parallel hyperplanes $n\cdot x=2\tau^{-2}$, which bound its central core $\tau^{-4}\B$, again a
$120$-cell.  The pieces are $2\,455\,201$ cells in $338$ classes, listed with exact coordinates in
the catalogue of \cite{SDR}; class $1$ is the core.  Our geometric search, run on this catalogue,
found $718$ cluster rules by an exhaustive search in exact coordinates over the scales $\tau^k$
($0\le k\le6$), $\sqrt5$, $\sqrt5\tau$ and $2$.  Here the linear system of
Theorem~\ref{thm:rules}, with the levels $\pm2\tau^2$, $\pm2\tau^{-2}$, is solved for the same
scales and for $\sqrt5\tau^2$, with no geometry beyond the Gram matrix of the normals
(\texttt{tools/zfoam4.py}, $32$ minutes in pure Python).

\paragraph{The module.}  The catalogue has $240$ hyperplanes: $120$ facets $n\cdot x=2\tau^2$ and
$120$ cuts $n\cdot x=2\tau^{-2}$, with $n$ twice a unit icosian, $|n|^2=4$.  The ratio of the two levels is the unit $\tau^4$, the
inflation that carries the central core onto the seed $\B=\tau^4\cdot$core (in three dimensions
the ratio was $\tau^3$).  The $\Zt$-span $\M_4$ of the $120$ normals has the Hermite basis
\[
(1,1,1,1),\quad(0,\,1-\tau,\,-\tau,\,-1),\quad(0,0,\,2\tau,\,2\tau-2),\quad(0,0,0,\,2-2\tau)
\]
(reduction over the norm-Euclidean ring $\Zt$), of determinant $4\tau-4$ and norm $-16$: the
index of $\M_4$ in $\Zt^4$ is $16$.  Two facts differ from three dimensions.  The $\Z$-span of the
$120$ normals is all of $\M_4$ and is $\tau$-invariant, so there is no analogue of the
distinction between $\Pm$ and $\M$.  And only $2\,593$ of the $3\,461$ vertices of the class
representatives lie in $\M_4$ (all $600$ vertices of the $120$-cell do); the cells of the star cut
have vertices in a finer module.  The argument for $t\in\M$ in Theorem~\ref{thm:rules} therefore
does not apply, and the program tests each translation: every translation of every rule below
lies in $\M_4$.

\paragraph{The group and the cells.}  The linear maps permuting the $120$ normals are found from
the images of four independent normals with the same Gram matrix: $14\,400$ elements, the full
symmetry group of the $120$-cell ($7\,200$ of determinant $+1$ and $7\,200$ of determinant $-1$).
Since they permute the generators of $\M_4$, their matrices in its basis have entries in $\Zt$.  The orbits of the $338$ class sign vectors
under the permutations of the $120$ cuts have the sizes recorded in the catalogue and give the
$2\,455\,201$ cells of the arrangement as $120$-bit integers; a union of cells is then a mask
test.

\paragraph{The rules.}  For each class and each scale $\lambda\in\{\tau^0,\dots,\tau^6\}\cup
\{\sqrt5,\sqrt5\tau,\sqrt5\tau^2\}\cup\{2\}$, four facets with independent normals and the
$4^4$ level assignments give the candidate translations; the surviving ones are those with every
facet on a level with the right orientation and the copy inside $\B$.  The volume identity is then
checked exactly for the cells selected by the sign constraints.  Unlike in
Result~\ref{res:rules}, single-cell solutions other than the identity are counted.

\begin{result}\label{res:4d}
The linear system has $853$ solutions $(A,\lambda,t)$ with $t\in\M_4$ for all of them, forming
$462$ pairs (class, scale) and $730$ triples (class, scale, composition).  By scale, the pairs
are: $\lambda=1$: $85$; $\tau$: $8$; $\tau^2$: $12$; $\tau^3$: $254$; $\tau^4$: $9$; $\tau^5$: $5$;
$\tau^6$: $65$; $\sqrt5$: $11$; $\sqrt5\tau$: $6$; $\sqrt5\tau^2$: $6$; $2$: $1$.  All $718$
triples of the geometric search are among them, with identical compositions; the $12$ others are
the $8$ triples at the scale $\sqrt5\tau^2$, which that search did not cover, and $4$ single-cell
translations at $\lambda=1$ (the classes $301\leftrightarrow304$ and $321\leftrightarrow322$ are
translates of each other).  $51$ of the $338$ classes have no rule at any of these scales.
\end{result}

The geometric search counted placements up to the symmetry of the class; the $853$ here are distinct
translations, which is why its counts are smaller for $183$ of the $452$ common pairs.  The
non-unit scales $\sqrt5\tau^j$ and $2$, absent in three dimensions (Result~\ref{res:rules}), do
occur: the facet-level condition $\lambda c+n\cdot t\in\{\pm2\tau^2,\pm2\tau^{-2}\}$ is a system
over $\Q(\sqrt5)$ and nothing forces $\lambda$ to be a unit.  In summary, the linear system of
Theorem~\ref{thm:rules} with the levels of the star cut reproduces the geometric search completely
and extends it to a scale that search had not covered; scales outside the list were not examined.
What does not carry over is the foam.  The substitution of Section~\ref{sec:subst}, the address
and the counting apply to the four-dimensional rules as written, but the cells of the $51$ classes
without a rule can never be subdivided, and they hold most of the volume.  In the figures below
each class with a rule uses a rule that contains a core when it has one, as in three dimensions.

\begin{figure}[tp]\centering
\fig[0.49\linewidth]{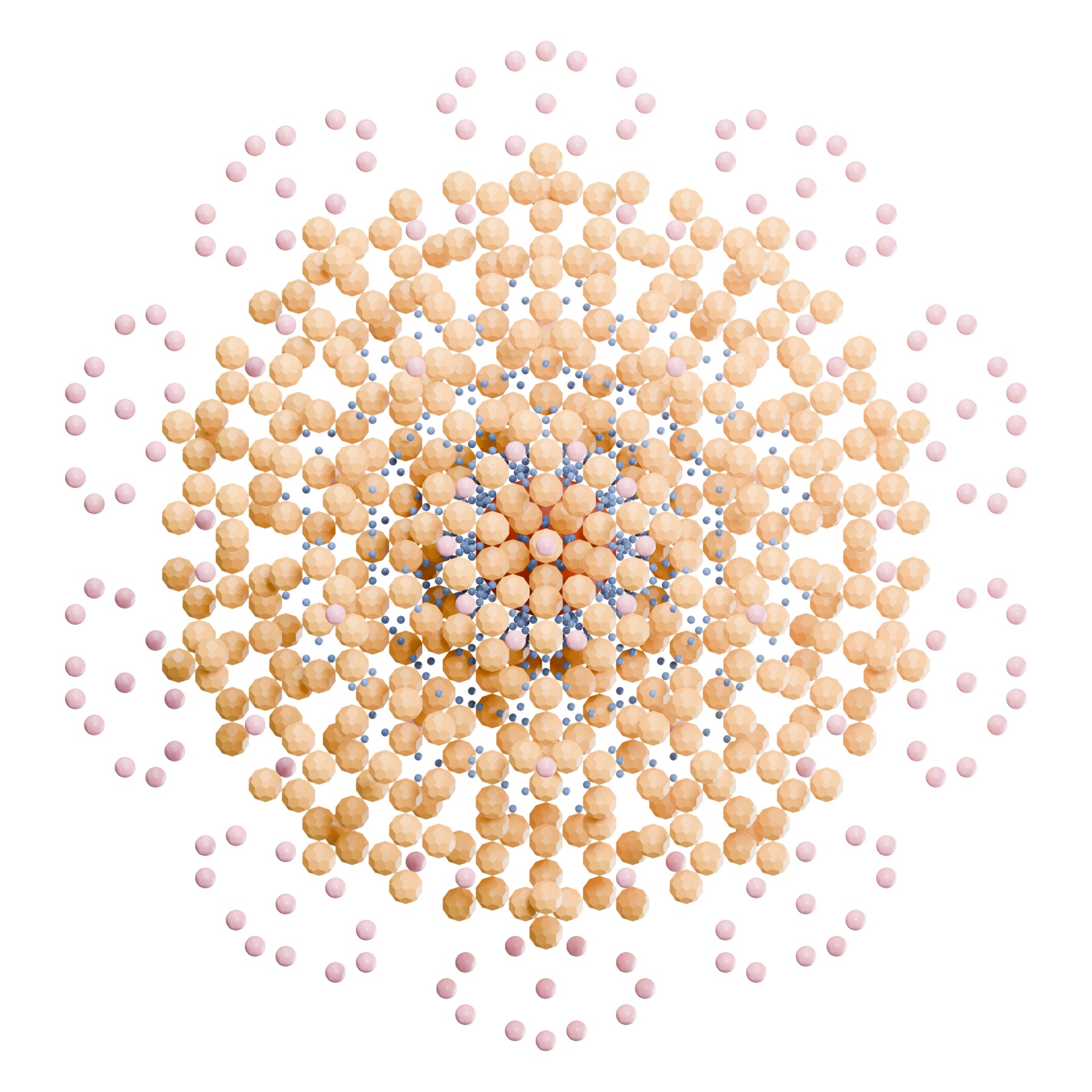}\hfill
\fig[0.49\linewidth]{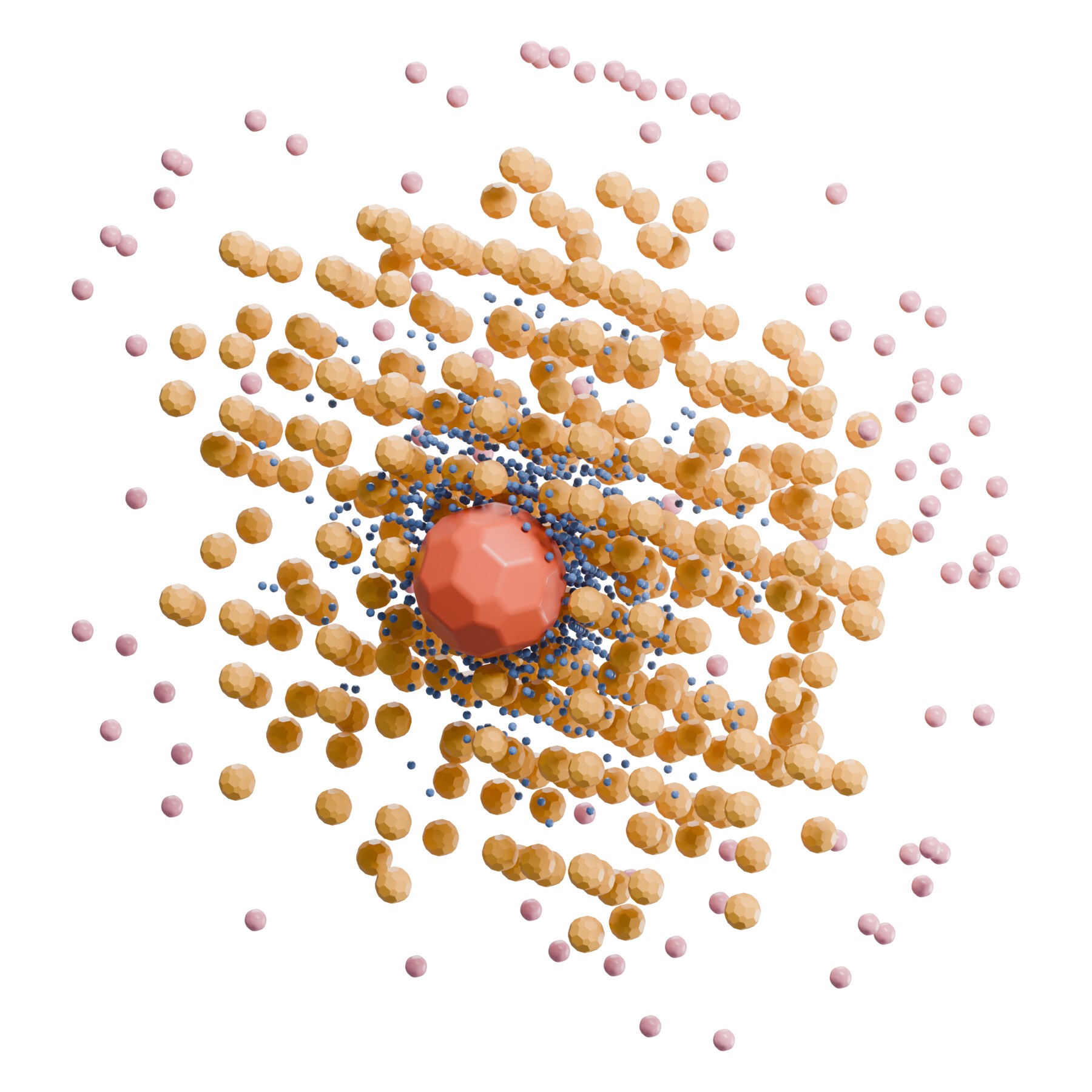}
\caption{The four-dimensional substitution after two steps: the \fdsecbub{} bubbles that meet
  the central section $n\cdot x=0$ of $\B=\tau^4\cdot$core, $n$ a facet normal.  Red: the core;
  gold: \fdbubthree{} bubbles at scale $\tau^{-3}$; pink: \fdbubfour{} at $\tau^{-4}$, on the outer
  shell just inside the facets of $B$; blue: \fdbubsix{} at $\tau^{-6}$.  Left: whole, seen along
  a five-fold axis of the section (the stabiliser of a facet normal in the symmetry group of the
  $120$-cell is the icosahedral group).  Right: cut in half, from the viewpoint of
  Figure~\ref{fig:foam4d}, with the core at the centre.  The bubbles fill $\fdvbub\,\%$ of this
  section; those of the three-dimensional Dodecafoam fill $56.5\,\%$ of space after two
  steps.}\label{fig:foam4dbub}
\end{figure}

\begin{figure}[p]\centering
\makebox[\linewidth][c]{\includegraphics[width=1.18\linewidth]{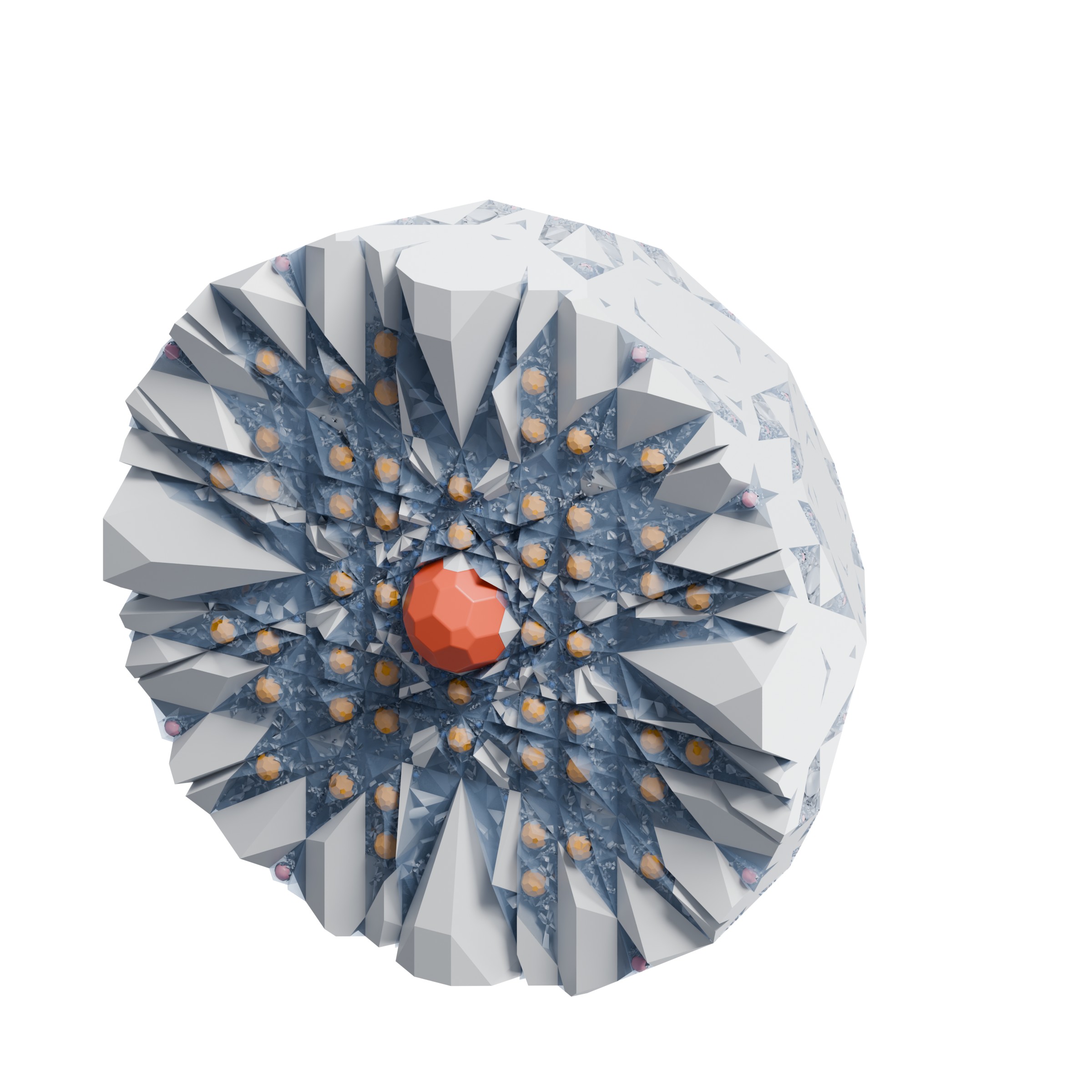}}
\caption{The whole of the central section of Figure~\ref{fig:foam4dbub}, cut in half.  Grey: the
  \fdsecfrozen{} frozen cells (ruleless classes; light grey from the first step, darker from the
  second), $\fdvfrozen\,\%$ of the section.  Translucent blue: the \fdsecactiveone{} active cells of
  the first step, which contain the \fdsecactive{} active tiles of the second step
  ($\fdvactive\,\%$) that are too small to draw.  In colour, the \fdsecbub{} bubbles
  ($\fdvbub\,\%$), coloured by scale as in Figure~\ref{fig:foam4dbub}.  Computed with the geometric
  substitution program on the \fdsectiles{} tiles that meet the hyperplane; the full tile set at this depth,
  $\fdfulltwo$ tiles, is never enumerated.}\label{fig:foam4d}
\end{figure}

\paragraph{A picture, and why there is no picture of the whole.}  Figures~\ref{fig:foam4dbub}
and~\ref{fig:foam4d} show the four-dimensional substitution after two steps.  They are
sections, not projections, and they are sections of a small part of the foam; both restrictions
are forced by the numbers.  The three-dimensional Dodecafoam has $165$, $2\,729$, $43\,933$ and
$708\,417$ tiles after one to four steps, a factor of about $16$ per step.  In four dimensions the
first step alone produces the $2\,455\,201$ cells of the star cut, the second $\fdfulltwo$ tiles
and the third $\fdfullthree$ (by the count matrix of the rules): factors of about $3\,000$ and
$7\,000$, because a cluster at scale $\lambda$ contains $\lambda^4$ times the cell density of the catalogue
cells, and the catalogue has $2\,455\,201$ cells where the dodecahedron has $165$.  The tiles
of the second step are too many to store, and the freezing does not help: the $51$ classes
without a rule stop refining, so most of the count is spent on ever smaller
active cells inside the $33\,\%$ of the volume that still moves.

What can be drawn is a three-dimensional section.  Take the hyperplane $n\cdot x=0$ through the
centre of $B$, $n$ the normal of a facet; it is the central section of \cite{SDR}, which
cuts the star cut into \fdsecpieces{} pieces.  A tile that does not meet the
hyperplane has no descendant that meets it, so the substitution can be iterated on the tiles that
straddle the hyperplane alone and the pruning is exact (\texttt{chok\_subst --section}, the geometric
substitution program with one filter added).  After two steps the section contains \fdsectiles{} tiles
(\fdsectime{} minutes): \fdsecbub{} bubbles, \fdsecfrozen{} frozen cells and \fdsecactive{} active
ones.  The active tiles are $92\,\%$ of the count and $\fdvactive\,\%$ of the volume of the section,
most of them at scales $\tau^{-3}$ to $\tau^{-6}$ of the catalogue cells, below the size of a pixel;
a third step in the section would have of the order of $10^{10}$ pieces.  So the picture is drawn at
depth $2$ with a convention: the frozen cells of both depths and the bubbles are drawn as they are,
and the active cells are drawn as their depth-$1$ parents, translucent---the union of the children
of a cell is the cell, so the solid shown is exactly the depth-$2$ section, with the parts that
would still change marked by the translucency.

The volume fractions of the section (three-dimensional volume) are, after one step, frozen
$\fdvfrozenstart\,\%$, active $\fdvactivestart\,\%$ and the central core $\fdvcore\,\%$; after
two steps, frozen $\fdvfrozen\,\%$ ($\fdvfrozenone\,\%$ in cells of the first step,
$\fdvfrozentwo\,\%$ in cells of the second), bubbles $\fdvbub\,\%$ and active $\fdvactive\,\%$.
The four-dimensional fractions are $65.5\,\%$, $1.6\,\%$ and $32.8\,\%$ (the central section has about twice the average share of bubbles).  The bubbles come in three sizes besides
the core: \fdbubthree{} at scale $\tau^{-3}$ ($\fdvthree\,\%$ of the section), \fdbubfour{} at
$\tau^{-4}$ ($\fdvfour\,\%$, all on the outer shell, just inside the facets of $B$) and
\fdbubsix{} at $\tau^{-6}$ ($\fdvsix\,\%$, the dust around the core).  Compare the three-dimensional
foam after two steps: $2\,729$ tiles, bubbles $56.5\,\%$ of the volume.  The foam has the full
symmetry of the $120$-cell, since the rule of a cell is the rule of its class transported by the
group element that places the cell, and the section has the icosahedral symmetry of the stabiliser
of $n$; Figure~\ref{fig:foam4dbub} looks down a five-fold axis.  What one sees of a bubble is the
polyhedron in which the hyperplane cuts a $120$-cell, with icosahedral symmetry for the central
one and smaller, less regular polyhedra for the others.

\section{Discussion}\label{sec:discussion}

\paragraph{Summary.}  Our geometric search found the Dodecafoam rules by geometry and verified
them in $\Zt$.  The reformulation shows that the rules are statements in $\Zt$: the arrangement
has the levels $1$ and $\tau^3$, whose ratio is a unit; inflation by a unit maps $\M$ onto itself; a rule is
a solution of $\lambda c+e\cdot t\in\{\pm1,\pm\tau^3\}$ inside $\B$; the group acts by matrices in
$\GL(3,\Zt)$ in the module's own basis; a tile is six integers, a group index and a class.
Returning to the question of the introduction: isometries are needed only as permutations of the
twelve normals and integer matrices; iteration is not needed for the question whether $x$ is a
bubble centre (the address) or for the production of a layer (the sieve); but there is no
cut-and-project formula.  Unlike the icosahedral model sets, such as the Ammann--Kramer--Neri
tiling \cite{BaakeGrimm12}, the Dodecafoam has bubbles of infinitely many sizes, those of one size
become sparse like $(\beta/\tau^3)^k$, and their window has measure zero.  The growth rate $\beta$
of the bubble counts and the decay rate $\rho$ of the uncovered volume are algebraic numbers
outside $\Zt$; $\beta$ is not in $\Q(\sqrt5)$.

\paragraph{Open points.}  (i) Whether $\dimH\Om=\log_\tau\beta$; this would follow from an open
set condition for the graph-directed system.  (ii) Whether $\Om_0$, $\Om_1$, $\Om_2$ have the same
limit points.  (iii) Whether the lift of $\M$ meets $\Om$ at points that are not bubbles; we
expect that it does.  (iv) The four-dimensional foam: the equivalence of
Theorem~\ref{thm:rules} holds in four dimensions (Section~\ref{sec:4d}), but $51$ classes have no
rule, and whether another choice of levels or of the cut gives a four-dimensional foam without
frozen classes is open.

\paragraph{Software and data.}  The programs are distributed with this article as ancillary files.
\texttt{tools/zfoam.py} (subcommands \texttt{check}, \texttt{rules}, \texttt{iterate},
\texttt{address}, \texttt{counts}, \texttt{sieve}, \texttt{window}, \texttt{test}, \texttt{bins})
covers Sections~\ref{sec:module}--\ref{sec:window}, together with \texttt{allscales.py}
(Result~\ref{res:rules}), \texttt{cone\_bound.py} (Proposition~\ref{prop:cone}) and
\texttt{window\_stats.py}.  \texttt{compare.py} produces the measurements of
Section~\ref{sec:compare}; \texttt{otherfoams.py}, \texttt{variant\_check.py}, \texttt{pair\_check.py}
and \texttt{foamlist.py} those of Section~\ref{sec:other}; \texttt{zfoam4.py},
\texttt{foam4d\_section.py} and \texttt{foam4d\_volumes.py} those of Section~\ref{sec:4d}, the
last two from the output of the geometric substitution program \texttt{chok\_subst}.  The logs of
all runs quoted in the text are in \texttt{out/}.

\end{document}